%% file: Revision1.tex
\documentclass[twocolumn, amsthm]{autart}

\usepackage{textcomp}
\usepackage{transparent}
\usepackage{stfloats}
\usepackage{url}
\usepackage{verbatim}
\usepackage{graphicx}
\def\BibTeX{{\rm B\kern-.05em{\sc i\kern-.025em b}\kern-.08em
    T\kern-.1667em\lower.7ex\hbox{E}\kern-.125emX}}
\usepackage{balance}

\usepackage{epsfig} 

\usepackage{epsfig} 
\usepackage{amsmath} 
\usepackage{amssymb}  
\usepackage{subfig}
\usepackage{algorithm}
\usepackage{multicol}
\usepackage{color}
\usepackage{threeparttable}
\usepackage[export]{adjustbox}
\usepackage{graphicx} 
\usepackage{caption,subcaption}
\usepackage{import}

\newtheorem{proposition}{Proposition}
\newtheorem{lemma}{Lemma}
\newtheorem{remark}{Remark}

\newtheorem{corollary}{Corollary}

\renewcommand{\vec}[1]{\ensuremath{#1}}

\usepackage{scalerel}

\newcommand\reallywidehat[1]{\arraycolsep=0pt\relax%
\begin{array}{c}
\stretchto{
  \scaleto{
    \scalerel*[\widthof{\ensuremath{#1}}]{\kern-.5pt\bigwedge\kern-.5pt}
    {\rule[-\textheight/2]{1ex}{\textheight}} 
  }{\textheight} %
}{0.5ex}\\           
#1\\                 
\rule{-1ex}{0ex}
\end{array}
}

\begin{document}
\begin{frontmatter}
\title{On Data-Driven Model Identification for Nonlinear Optimal Control}

\author[AS]{Aayushman Sharma}\ead{aayushmansharma@tamu.edu},
\author[AS]{Suman Chakravorty}\ead{schakrav@tamu.edu}

\address[AS]{Department of Aerospace Engineering, Texas A\&M University, College Station, TX 77843 USA.}


\begin{abstract}
In this paper, we study the use of nonlinear model identification techniques for the optimal control of nonlinear systems, also known as model-based Reinforcement Learning. We show that the nonlinear model identification problem is equivalent to estimating the generalized moments of an underlying sampling distribution and is bound to suffer from ill-conditioning and variance when approximating a system to high order and over a large domain, requiring samples combinatorial-exponential in the order of the approximation and domain size: a ``Curse of Variance and Ill-Conditioning (COVIC)" that shows up even in very low dimensional problems, quite apart from the usual ``Curse of Dimensionality". We show that the iterative identification of ``local" linear time varying (LTV) models around the current estimate of the optimal trajectory, coupled with a suitable optimal control algorithm such as iterative LQR (ILQR), alleviates these issues and is sufficient to locally accurately solve the underlying optimal control problem.

\end{abstract} 
\end{frontmatter}
\input{Introduction.tex} 
\input{Theory_revised}
\input{Simulation_new}

\input{Bibliography.tex}

\end{document}

%% file: Introduction.tex
\section{Introduction}

\begin{figure}[!htb]


      \subfloat[15-Link Swimmer Robot]{\includegraphics[width=0.5\linewidth]{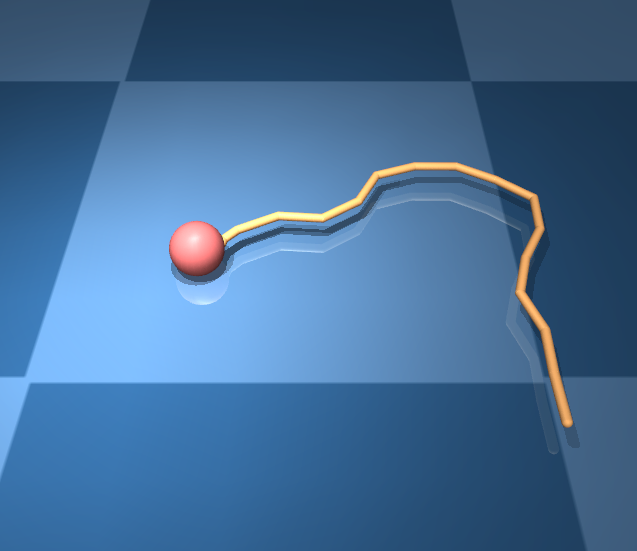}}    
      \subfloat[Fish Robot]{\includegraphics[width=0.5\linewidth]{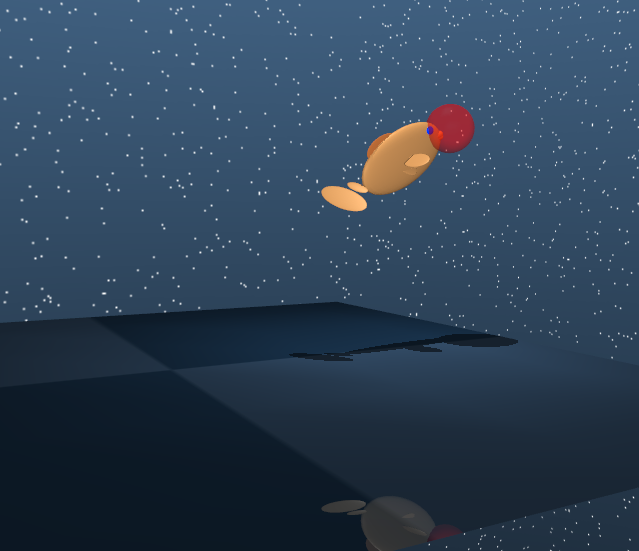}}
      

\caption{{In many scenarios, the dynamics of the physical systems are generally analytically unknown, and we only can access an input-output computational `black box'. In the swimmer examples above, the fluid-structure interactions lead to analytically unknown equations.}}
\vspace{-0.4cm}
\label{fig:unknowndynamics}
\end{figure}

Many physical systems are analytically intractable (Fig.~\ref{fig:unknowndynamics}). 
Thus, in order to devise strategies for optimal control of such systems, there is a need for accurate identification of the dynamics of these systems. This is especially attractive in scenarios where utilization of a `black-box' computational model is expensive, and thus, a simplified `surrogate' model identified from the original dynamics can be used to replace the actual system. 
In this paper, we study the use of data-driven surrogate modeling for nonlinear optimal control problems. We assume the availability of full state information in this work.

\begin{figure}[!htpb]
    \centering
    \includegraphics[width=1\linewidth]{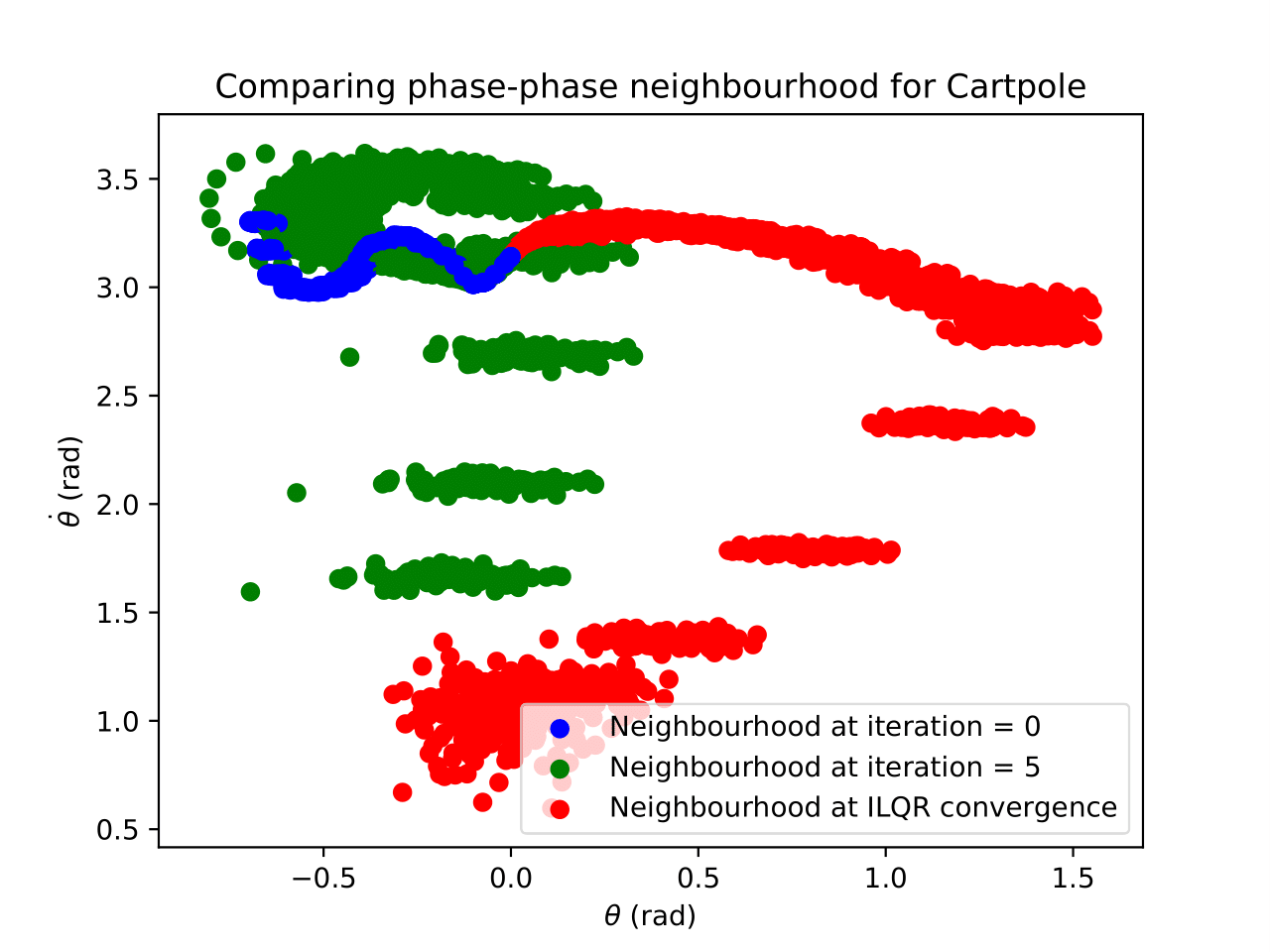}
    \caption{For optimal control problems, the neighbourhood of the optimal trajectory may be significantly different from the neighbourhood of the initial guess. Thus, a surrogate model trained in the vicinity of a random initialization may not generalize to the optimal trajectory. Hence, there is no well-defined domain of interest in this class of problems unless it is assumed to be ``large enough" to encompass the true optimum.}
    \label{fig:placeholder}
\end{figure}

The research topic of interest in this context is nonlinear system identification. Nonlinear system identification is a well-researched topic, and various approaches such as: block-oriented structures \cite{giriblock}\cite{schoukens2015parametric}, nonlinear state-space representations \cite{paduart2010identification}\cite{schon2011system}, nonlinear autoregressive exogeneous (NARX), nonlinear autoregressive moving average with exogenous input (NARMAX) models \cite{billings2013nonlinear}, and piecewise linear models \cite{mattsson2016recursive}, have been proposed to address the problem. 
In addition, neural-net based approximators such as deep neural networks, RNNs and LSTMs \cite{narendra1989identification}\cite{dinh2010dynamic}\cite{ogunmolu2016} have been used to identify models.
In spite of this research, it is well understood that nonlinear identification problems are very complex and difficult to solve, and require addressing three fundamental issues: 1) experimental design/ data gathering, 2) the nonlinear model structure used, and 3) the model estimation approach, as succinctly noted in reference \cite{schoukens2019nonlinear}. In this paper, we address items 1-3 above in the context of nonlinear optimal control. A good nonlinear identification approach requires careful data gathering over the ``domain of interest" \cite{schoukens2019nonlinear}. However, general optimal control problems do not have a pre-defined domain of interest (see Fig.~\ref{fig:placeholder}). Thus, any identification approach needs to search over a suitably large region that can be assumed to contain the optimal solution resulting in the need for prohibitive amounts of data that is exponential-combinatorial in the model order and size of the domain, regardless of the model structure and even in the presence of full-state measurements. Furthermore, this results in a highly ill-conditioned model estimation problem leading to inaccurate solutions. This problem shows up even in very low-dimensional problems and is a ``Curse of Variance and Ill-Conditioning (COVIC)" quite apart from the usual ``Curse of Dimensionality" \cite{bellman1957}.\\

The basic problem of controlling nonlinear systems is very challenging. Despite a rich history of control techniques for nonlinear dynamical systems, the area of data-driven optimal control design is still a topic of very active research called Reinforcement Learning (RL) \cite{lillicrap2015continuous}\cite{schulman2015trust}\cite{schulman2017proximal}. 
In the case of model-based reinforcement learning, many leading techniques seek to identify and build a unified, 'global' surrogate model to predict the dynamics across a large portion of the state space. Algorithms such as PETS \cite{pets} and MBPO \cite{mbpo}, for instance, learn an ensemble of neural networks to serve as a comprehensive dynamics model for planning and data augmentation. Similarly, architectures like PILCO \cite{pilco} and MC-PILCO \cite{mcpilco} use a rollouts-based approach to identify a `global' Gaussian Process-based model that is capable of prediction in the entire state space. More advanced `world models', such as the Dreamer architecture \cite{dreamer,dreamerv3}, take this to the next level by building a predictive model from a large experience replay buffer of past samples, aiming to generalize and enable agents to learn complex behaviors entirely within this learned ``dream."
However, this pursuit of a single, global model is fundamentally ill-posed owing to the issues of ill-conditioning and variance that we identify in this paper. Any finite replay buffer is merely a sparse sampling of the true state space, and the resulting models are not truly global since the samples required for a good solution are combinatorial-exponential in the size of the approximation domain and the order of the approximation. Their success hinges on interpolating within the training data distribution, but when forced to interpolate to new (unseen) states, yield unreliable predictions that destabilize the control optimization process. The challenges inherent in these  `global' methods motivate the need for a more robust and reliable approach.

In previous work \cite{wang2021search}\cite{D2C_ASME}, we proposed a model-based RL algorithm called Decoupled Data-based Control (D2C) for optimal nonlinear control, which identifies successive LTV models combined with the so-called ILQR approach \cite{ilqg2} to yield a globally optimum local feedback law and is far more efficient, reliable and accurate when compared to state of the art RL techniques \cite{D2C_ASME}. Thus, in the context of data-driven nonlinear optimal control, we show that such a recursive local LTV scheme suffices to find accurate solutions to the optimal control problem as opposed to the ``global" approach. \\

\textbf{Contributions.} In this paper, we show that: 1) a nonlinear model identification problem is equivalent to a generalized moment estimation problem of an underlying sampling distribution, and when approximating a system to high order over a large domain, the resulting model estimation problem is bound to suffer from ill-conditioning as well as high variance, that scales in a combinatorial-exponential fashion in the size of the domain and the order of the approximation, regardless of model structure; 
2) extensive empirical evidence is provided to support the theoretical claim showing that the Gram matrix $\mathcal{G}$ is highly ill-conditioned while there is high variance in the forcing function $F$ in the Normal equation $\mathcal{G}C = F$ underlying the model estimation problem leading to unreliable identification of the model parameters $C$; and  
3) in the context of the target problem of solving nonlinear optimal control problems, the use of a recursive ``local" LTV identification approach, as opposed to ``global" nonlinear approach as above, is sufficient 
to address the issues of ill conditioning and variance, assuring accurate and provably optimum local answers.
This work is a significantly expanded version of our conference paper \cite{conf_paper}. It builds upon the prior work by developing a generalized framework to formally characterize the approximation error through the ill-conditioning and variance, and provides significantly more empirical evidence for the theory. Furthermore, we introduce an analysis for the multidimensional case to demonstrate how the ill-conditioning and variance scale with the state dimension.

The rest of the paper is organized as follows: In Section \ref{probform}, we introduce the nonlinear optimal control problem, and the nonlinear model-identification problem. We discuss the issues of model convergence, ill-conditioning and variance of these techniques, regardless of architecture, and provide theoretical justifications for the same, in Sections \ref{truesec3} and \ref{sec_var_illcon}. In Section 5, we examine the implications of the ill-conditioning and variance inherent to model estimation on the optimal control problem and show that a recursive local LTV model estimation approach, allied with a local optimal control approach such as ILQR, suffices to alleviate the COVIC issue and accurately solve for the optimal control. In Section \ref{sec4}, we provide empirical evidence to support the proposed theory, wherein we benchmark the ``global" surrogate-models learned by the algorithms, and their performance on optimal control law computation as compared to a ``local" linear time-varying approach. In section 7, we draw conclusions on this work and discuss its ramifications.

%% file: Theory_revised.tex
\section{Preliminaries}
\label{probform}
Consider the discrete-time nonlinear dynamical system:
\begin{equation}
    x_{t+1}=f(x_t, u_t),
\label{eq:dyn}
\end{equation}
where $x_t\in \mathbb{R}^{n_x}$ and $u_t\in \mathbb{R}^{n_u}$ correspond to the state and control vectors at time $t$. The optimal control problem is to find the optimal control policy $\pi^0=\{\pi^0_0, \pi^0_1 ... \pi^0_{T-1}\}$, that minimizes the cumulative cost:
\begin{align*}
    & \min_\pi J^\pi(x)=\sum_{t = 0}^{T-1} c_t(x_{t}, u_{t}) + c_T(x_T),
\end{align*}
subject to the dynamics above,
given some $x_0 = x$, and where $u_t=\pi_t^0(x_t)$, $c_t(\cdot)$ is the instantaneous cost function and $c_T(\cdot)$ is the terminal cost. We assume the incremental cost to be quadratic in control, such that $c_t(x_t,u_t)=l_t(x_t)+\frac{1}{2}u_t^T R u_t$.

In this work, we assume that we only have access to the input-output model of the dynamical system, wherein given the inputs $x$ and $u$, a black-box (computational) model generates the next state $x'=f(x,u)$. Thus, we seek to find a surrogate model for the dynamics from the system input-output data, which can then be utilized to solve the nonlinear optimal control problem posed above.

\subsection{Model Structure} 
We consider three alternatives for the surrogate modeling:
\begin{enumerate}
\item a ``global" linear in parameters approximation architecture with nonlinear basis functions ($f (x,u) \approx \sum_i c_i \psi_i(x,u)$), 
\item a ``global" nonlinear in parameters architecture, such as a (deep) neural network, ($f(x,u) \approx h(x,u;\theta)$), and 
\item a ``local" linear time varying (LTV) architecture ($\delta x_t = A_{t-1}\delta x_{t-1} + B_{t-1}\delta u_{t-1}$), linearized around a nominal trajectory $(\bar{x}_t,\bar{u}_t)$.
\end{enumerate}
The first two architectures estimate a ``global" model of the underlying dynamics while the third only finds a ``local perturbation" model of the dynamics. Since we are considering an optimal control problem, the ``local" model is an LTV approximation of the perturbation dynamics around the current estimate $(\bar{x}_t,\bar{u}_t)$ of the optimal trajectory, i.e., $\delta x_t = x_t - \bar{x}_t, \delta u_t = u_t -\bar{u}_t$. Moreover, these LTV models need to be recursively updated as the optimization proceeds (see Fig.~\ref{fig:optcontrol}(a)-(b)). In the case of the global approximators, in theory, such a model update is not required, since the approximate model should have global validity given enough training data.  

For solving the optimal control problem using LTV models, we utilize the Iterative Linear Quadratic Regulator (ILQR) algorithm, which can be guaranteed to converge to the global optimum \cite{wang2021search}. A brief description of the ILQR Algorithm is given in the appendix.

\subsection{Data Generation and the Model Estimation Problem}

Let $f, g:\mathbb{R}\rightarrow \mathbb{R}$ be an infinitely smooth function, and $\pi(\cdot)$ be a p.d.f./sampling distribution on $\mathbb{R}$ from which training/ identification data is sampled. The inner product w.r.t. the distribution $\pi(\cdot)$ is defined as 
\begin{equation*}
   \langle f,g \rangle = \int_{-\infty}^{\infty} f(x)g(x)\pi(x)dx. 
\end{equation*}

We are interested in approximating the function $f(\cdot)$ using a suitable architecture. For simplicity, we consider a scalar, control-free dynamical model $f(\cdot)$. Nonetheless, the results are generalized in a relatively straightforward fashion, while the issues inherent in the estimation problem come into sharp focus in the one dimensional setting.

\textbf{Linear.} Given a set of basis functions $\{\psi_0 ... \psi_M\}$, we want to approximate $f$ in the subspace spanned by $\Psi = \{\psi_0...\psi_M\}$, i.e., we want to find 
\begin{equation}
    \hat{f}=\sum_{i=0}^{M} \hat{c}_i \psi_i,
\end{equation}
such that $||f-\hat{f}||$ is minimized, where $||f||=\sqrt{\langle f,f \rangle}$ is the norm induced by the $\pi(.)$ weighted inner product above. Further, we want to find the number of samples $R$ required such that an $\varepsilon-$close approximation can be made to the true solution given we only have access to samples $(x^i,f(x^i))$ from the underlying distribution.

\textbf{Nonlinear.} Given a nonlinear approximation architecture $h(x,\theta)$, such as a (deep) neural network, we want to find the best approximation of $f$ in the approximation architecture, i.e. find $f^* = h(x,\theta^*)$ such that $||f-f^*||^2$ is minimized over all parameters $\theta \in \Theta \subseteq \Re^n$. Further, we want to approximate the number of samples to make an $\varepsilon-$close approximation based on only the samples $(x^i,f(x^i))$ from the underlying distribution.

\subsection{Heuristic Overview}
In the following, we give a brief overview of the developments in Sections \ref{truesec3}-\ref{section5}. 
In Sections \ref{truesec3} and \ref{sec_var_illcon}, we shall concentrate purely on
the model estimation problem posed above and its properties while we will examine the implications of this analysis on the nonlinear optimal control problem in Section \ref{section5}. In particular, suppose that the size of the sampling distribution $\pi(.)$ is parametrized by a scalar variable $\sigma$.

In Section \ref{truesec3}, we will show that the linear in parameters architecture converges to the ``correct" model given sufficient data from the aformentioned sampling distribution. 
We also give sample complexity bounds on the number of samples $R$ required for a given level of accuracy with some desired confidence. 
The nonlinear in parameters model analysis follows a similar construction in a recursive form and we give similar sample complexity bounds for the R-sample calculation of the gradient. 

In Section \ref{sec_var_illcon}, we show that the number of samples $R$ required for an accurate solution scales combinatorially-exponentially in the size of the sampling distribution $\sigma$ and the implicit order of the approximation required in the sampled domain, say $N$ (defined in Section 4), owing to high variance and severe ill-conditioning in the underlying least squares model estimation problem. 
A similar sample complexity is shown to hold for every gradient descent step in the nonlinear model estimation problem. In practice, this ``Curse of Variance and Ill-Conditioning (COVIC)" makes the estimation of a global model infeasible.

In Section \ref{section5}, we return to the nonlinear optimal control problem. The lack of an accurate ``global" surrogate model implies that we cannot solve the problem utilizing the aforementioned one-shot global data-driven models. In lieu, utilizing the ILQR approach to solving optimal control problems, it naturally follows that ``limiting" the model estimation to a recursive ``local" LTV identification suffices to accurately solve the optimal control problem. Furthermore, we show that this implies that any data-based approach to solving nonlinear optimal control problems will have to rely on recursive local model identification, but which need not necessarily be LTV models.


\section{Convergence in the Model Estimation Problem}\label{truesec3}
In this Section, we characterize the solution to the model estimation problem posed above in the linear as well as nonlinear case and establish convergence given sufficient data.
\subsection{Linear Case}
We can get the answer to the linear question from the Projection Theorem in the Hilbert space $L_2(p)$ \cite{Courant-Hilbert}. 

The optimal $\hat{f}$ is such that $(f-\hat{f})\perp \psi_i \forall i=1,2...M$, owing to the Projection Theorem. Substituting $\hat{f} =\sum_{i=0}^M \hat{c}_i \psi_i$ leads to the Normal equations:
\begin{align}
    \mathcal{G} \hat{C}_M &= F, \label{N}\\
    \mathcal{G}_{ij} &= \gamma_{ij} = \langle \psi_i, \psi_j \rangle
    = \mathbb{E}[\psi_i\psi_j] \nonumber \\
    f_i &= \langle f, \psi_i \rangle 
         = \mathbb{E}[f \, \psi_i],  \nonumber
\end{align}
where the matrix $\mathcal{G} = [\mathcal{G}_{ij}]$, and the vector $F = [f_i]$, $\hat{C}_M= [\hat{c}_0, \cdots \, \hat{c}_M]'$, and $M$ is the order of the approximation.

Next, we assume that $p(x)$ is specified by an exploration parameter $\sigma$. For instance:
\begin{align}
    \pi(x) &\sim \mathcal{N}(x; 0, \sigma^2) \quad \text{(Gaussian)}, \nonumber \\
    \pi(x) &\sim U(x; -\sigma, \sigma) \quad \text{(Uniform)}. \nonumber
\end{align}
Thus, the functions $f_i(.)$ are a function of the parameter $\sigma$.

Next, let us assume that $f$ is exactly approximated in the basis set given an infinite number of basis functions, and define:
\begin{align}
    f(x) &= \sum_{k=0}^\infty c_k^* \, \psi_k(x), \nonumber \\
    f_i(\sigma) &= \sum_{k=0}^\infty c_k^* \, \mathbb{E}[\psi_i \, \psi_k]
               \equiv \sum_{k=0}^\infty c_k^* \, \gamma_{i,k}(\sigma), \nonumber \\
    f_{i,M}(\sigma) &= \sum_{k \le M} c_k^* \, \gamma_{i,k}(\sigma), \nonumber\\
    \delta f_{i,M}(\sigma) &= \sum_{k > M} c_k^*
    \gamma_{i,k}(\sigma). \label{eq:delta}
\end{align}
Note the distinction between the coefficients $\hat{c}_k$ and $c_k^*$ and that $\hat{c}_k \rightarrow c_k^*$ as $M \rightarrow \infty$.\\

Next, we digress a little and have the following supporting development.
Let $g(\sigma)$ be an infinitely smooth function of the variable $\sigma \ge 0$, and further let  
\begin{align}
    g(\sigma) &= \sum_{k=0}^\infty g_k \, \zeta_k(\sigma). \nonumber
\end{align}

Define: 
\begin{align*}
    g_M(\sigma) = \sum_{k \le M} g_k \, \zeta_k(\sigma),\;
    \delta g_M(\sigma) = \sum_{k > M} g_k \, \zeta_k(\sigma). 
\end{align*}

\textbf{Definition:} The infinitely smooth function $g(\sigma)$ is said to have the \emph{Uniform Approximation (UA)} property with respect to the series of infinitely smooth functions $\{\zeta_k(\sigma)\}$ if  given any $M> 0$, and any $\varepsilon >0$, there exists $\sigma_M^{\varepsilon} > 0$ such that: $|\delta g_M(\sigma)| \leq \varepsilon |g_M(\sigma)|$, for any $\sigma \leq \sigma_M^{\varepsilon}$.

The above definition models the intuition that for an approximation scheme to be applicable, any $M^{th}$ order expansion of a function is sufficiently good in some suitable bounded region.
For instance, the UA is satisfied if the functions $\zeta_k(\sigma)$ are polynomials of degree $k$ and if the function $g(\sigma)$ has an infinite Taylor series expansion.\\



\begin{lemma}  \label{e-lemma}
Let the functions $f_i(\sigma)$ have the UA property with respect to the sequence of infinitely smooth functions $\{\gamma_{i,k}(\sigma)\}$ for any $i$.  
Then, there exists a $\bar{\sigma}$ such that for all $\sigma \le \bar{\sigma}$, 
\begin{align}
    |\delta f_{i,M}(\sigma)| \le \varepsilon \, |f_{i,M}(\sigma)|, \quad \forall i,
\end{align}
where the functions $f_{i,M}(\sigma)$ and $\delta f_{i,M}(\sigma)$ above have been defined in \eqref{eq:delta}.
\end{lemma}
\begin{proof}  
Since $f_i(\sigma)$ has the UA property, for each $i$ there exists $\sigma_{i,M}^{\varepsilon}$  
such that for all $\sigma \le \sigma_{i,M}^{\varepsilon}$,  
$
    |\delta f_{i,M}(\sigma)| \le \varepsilon \, |f_{i,M}(\sigma)|.
$
Take  
\begin{align}
    \bar{\sigma} = \min_i \sigma_{i,M}^{\varepsilon}, \nonumber
\end{align}
then it follows that  for any $i$,
\begin{align}
    |\delta f_{i,M}(\sigma)| \le \varepsilon \, |f_{i,M}(\sigma)|, \quad \sigma \le \bar{\sigma}. \nonumber
\end{align}
\end{proof}
Let  
\begin{align}
    C_M^* &= 
    \begin{bmatrix}
        c_0^* \\ \vdots \\ c_M^*
    \end{bmatrix},
    &
    \hat{C}_M &=
    \begin{bmatrix}
        \hat{c}_0 \\ \vdots \\ \hat{c}_M
    \end{bmatrix},
\end{align}
where $\hat{c}_k$ is define in \eqref{N} and $c_k^*$ is defined in \eqref{eq:delta}.
\begin{proposition}  \label{ME_acc}
Given the sampling parameter $\sigma \le \bar{\sigma}$,
\begin{align}
    \| C_M^* - \hat{C}_M \| &\le \kappa(\mathcal{G}) \, \| C_M^* \| \, \varepsilon,
\end{align}
where  
$
    \kappa(\mathcal{G}) = \| \mathcal{G} \| \, \| \mathcal{G}^{-1} \|
$
is the condition number of the Gram matrix, and $\| \cdot \|$ is the matrix norm induced by the vector norm $||.||$.
\end{proposition}

\begin{proof}  
Consider the model estimation equation:
\begin{align}
    \mathcal{G} \hat{C}_M &= F(\sigma) \equiv
    \begin{bmatrix}
        f_0(\sigma) \\
        \vdots \\
        f_M(\sigma)
    \end{bmatrix}, \nonumber
\end{align}
However,
$
    f_i(\sigma) = f_{i,M}(\sigma) + \delta f_{i,M}(\sigma).
$

Note that:
\begin{align}
    F_M(\sigma) &\equiv
    \begin{bmatrix}
        f_{0,M}(\sigma) \\
        \vdots \\
        f_{M,M}(\sigma)
    \end{bmatrix}
    =
    \mathcal{G} C_M^*. \nonumber
\end{align}

Hence, defining:
\begin{align*}
    \Delta F_M(\sigma) = 
    \begin{bmatrix}
    \delta f_{0,M} (\sigma)\\
    \vdots\\
    \delta f_{M,M} (\sigma)
    \end{bmatrix}
\end{align*}
it follows that:
\begin{align}
    \mathcal{G} \, \hat{C}_M &= \mathcal{G} \, C_M^* + \Delta F_M, \\
    \hat{C}_M - C_M^* &= \mathcal{G}^{-1} \Delta F_M. \nonumber
\end{align}

Therefore:
\begin{align}
    \| \hat{C}_M - C_M^* \| &\le \| \mathcal{G}^{-1} \| \, \| \Delta F_M \|, \nonumber \\
    &\le \| \mathcal{G}^{-1} \| \, \varepsilon \, \| {F}_M \| 
    \quad \text{(by Lemma 1)}, \nonumber \\
    &\le \| \mathcal{G}^{-1} \| \, \| \mathcal{G} \| \, \| C_M^* \| \, \varepsilon 
    = \kappa(\mathcal{G}) \, \| C_M^* \| \, \varepsilon. \nonumber
\end{align}
\end{proof}
\begin{remark}
    In order to make an $\varepsilon$-close approximation of a function to $M^{th}$ order by sampling, we need a certain exploration $\bar{\sigma}$ that is dependent on all the moments $f_i(\sigma) = \mathbb{E}[f\psi_i]$, $i=0,1,2,\cdots M$. Furthermore, the accuracy of the approximation is critically dependent on the conditioning $\kappa(\mathcal{G})$ of the Gram matrix.
\end{remark}

Next, we assume that we have access only to the data-points $(x^k, f(x^k)), k=1,2...R$, where the samples are drawn from the sampling distribution, $x^k\sim p(\cdot)$, independently. 
Given the sample set $(x^k, f(x^k))$, we approximate:
\begin{equation}
\label{eq: Leastsquares}
    f(x^k) \approx \sum_{j=1}^M \hat{c}_j \psi_j (x^k)~ \forall~ k=1,2...R.
\end{equation}
The least squares solution to Eq.~\ref{eq: Leastsquares} is given by:
\begin{equation}
\label{eq: Leastsquares solution}
    \hat{C}^R_M = (\underbrace{{\Psi}^{R'}{\Psi}^R}_{\mathcal{G}^R})^{-1} {\Psi}^{R'} F^R,
\end{equation}

 where, 
 \begin{align*}
     ({\Psi}^{R'}{\Psi}^R)_{ij}=\sum_{k=1}^R \psi_i(x^k)\psi_j(x^k) \equiv R\gamma_{ij}^R, \\
   ({\Psi}^{R'}F^R)_i = \sum_{k=1}^R \psi_i(x^k)f(x^k) \equiv Rf_i^R,
 \end{align*}
 where $(.)_i, (.)_{ij}$ represent the $i^{th}$ element of a vector and the $ij^{th}$ element of a matrix respectively.
   
 Owing to the law of large numbers \cite{gnedenko}, we have: 
\begin{equation*}
    \frac{1}{R} \sum_{k=1}^R \psi_i(x^k)\psi_j(x^k)\rightarrow \int \psi_i(x)\psi_j(x)\pi(x)dx=\underbrace{\langle \psi_i,\psi_j\rangle}_{\gamma_{ij}},
\end{equation*}
and
\begin{equation*}
    \frac{1}{R} \sum_{k=1}^R \psi_i(x^k)f(x^k)\rightarrow \int \psi_i(x)f(x)\pi(x)dx=\underbrace{\langle \psi_i,f\rangle}_{f_i} ,
\end{equation*}
as $R\rightarrow\infty$. This implies that 
\begin{align*}
    \hat{C}^R_M&=({\Psi}^{R'} {\Psi}^R)^{-1}({\Psi}^{R'}F^R)\\
    &=(\frac{1}{R}{\Psi}^{R'} {\Psi}^R)^{-1}(\frac{1}{R}{\Psi}^{R'}F^R)\\
    &\rightarrow(\mathcal{G})^{-1}F = \hat{C}_M \text{ as } R\rightarrow\infty.
\end{align*}
Therefore, the data-based solution converges to the true solution of the Normal Equation, as $R\rightarrow \infty$.\\
Albeit the Least-Squares estimate converges to the true solution of the Normal equation in the limit, we would like to find an estimate of the samples $R$ required for some given error level $\epsilon$. In this regard, we first present the following facts.

Let $q(x)$ be a function of the random variable $x$. Then $\Bar{q}\equiv E[q(x)]=\frac{1}{R}\sum_{j=1}^{R}q(x^j)$ as $R\rightarrow\infty$. Furthermore, if $R$ is large enough, due to the Central Limit Theorem \cite{gnedenko}, $\Bar{q}^R \equiv \frac{1}{R}\sum_{k=1}^{R}q(x^k) = N(\Bar{q}, \frac{\sigma_q^2}{R})$,
where $\sigma_q^2= Var[q(x)]$, and $N(\mu,\sigma^2)$ is a normal random variable with mean $\mu$ and variance $\sigma^2$. Next, the r.v. $\frac{\bar{q}^R - \bar{q}}{\sigma_q/ \sqrt{R}} \sim N(0,1)$. Thus:
\begin{align}
    Prob(|\bar{q}^R - \bar{q}|> \epsilon) = 1- \Phi(\frac{\epsilon\sqrt{R}}{\sigma_q^2}),
\end{align}
where $\Phi(a)$ (termed the erf) represents the probability that a standard Gaussian r.v. z, i.e., $z\sim N(0,1)$, takes values $|z|<a$.\\
\\ 
Recall that $f_i = \langle \psi_i, f \rangle$, and thus:
\begin{align*}
    f_i^R &\equiv \frac{1}{R}\sum_{k=1}^R \psi_i(x^k)f(x^k)
    \sim N(f_i,\frac{\nu_i^2}{R}), \\
    \nu_i^2 &=Var[\psi_i(x)f(x)].
\end{align*}

Further, recall that the elements of the Gram matrix are $\gamma_{ij}=\langle \psi_i,\psi_j \rangle$, and that $\frac{1}{R}\underline{\Psi}^{R'}\underline{\Psi}^R \rightarrow \mathcal{G}$ as $R\rightarrow \infty$. Thus, 
\begin{align*}
    \gamma_{ij}^R & \equiv \frac{1}{R}\sum_{k=1}^R \psi_i(x^k)\psi_j(x^k) \sim N(\gamma_{ij}, \frac{\nu_{ij}^2}{R}),\\
    \nu_{ij}^2 &= Var[\psi_i(x)\psi_j(x)].
\end{align*}
Then, it follows from the development above that:
\begin{align*}
    Prob(|f_i^R-f_i|>\epsilon)= 1-\Phi(\frac{\epsilon \sqrt{R}}{\nu_i}),\\
    Prob(|\gamma_{ij}^R-\gamma_{ij}|>\epsilon)= 1-\Phi(\frac{\epsilon \sqrt{R}}{\nu_{ij}}).
\end{align*}

The above development may now be summarized in the following result.\\

\begin{proposition}\label{prop2}
The solution of the data-based Least-Squares problem (Eq.~\ref{eq: Leastsquares solution}) $\hat{C}^R_M\rightarrow \hat{C}_M$, the solution of the Normal Equation \eqref{N}, as $R\rightarrow\infty$, in the mean square sense. Furthermore, the R-sample estimates $f_i^R \sim N(f_i, \frac{\nu_i^2}{R})$, $\gamma_{ij}^R\sim N(\gamma_{ij}, \frac{\nu_{ij}^2}{R})$, and $ Prob(|f_i^R-f_i|>\epsilon)= 1-\Phi(\frac{\epsilon \sqrt{R}}{\nu_i})$, $ Prob(|\gamma_{ij}^R-\gamma_{ij}|>\epsilon)= 1-\Phi(\frac{\epsilon \sqrt{R}}{\nu_{ij}})$.
\end{proposition}
Suppose, now we want a failure probability of at most $\delta$ that each $|\gamma_{ij}^R - \gamma_{ij}|> \epsilon,\, |f_i^R - f_i| > \epsilon$. Then, the number of samples required for the $f_i^R$ to reach this confidence level is:
\begin{align*}
    R_i = \nu_i^2 \frac{p^2(\delta)}{\epsilon^2},
\end{align*}
and the number of samples required for the $\gamma_{ij}^R$ is:
\begin{align*}
    R_{ij} = \nu_{ij}^2 \frac{p^2(\delta)}{\epsilon^2},
\end{align*}
where $p(\delta) = \Phi^{-1}(1- \delta)$. Let:
\begin{align}\label{sample_C}
    R = \max\left( \max_i R_i, \max_{i,j} R_{ij}\right).
\end{align}
Let $\mathcal{G}^R = \mathcal{G} + \delta \mathcal{G}^R$ and $F^R = F+ \delta F^R$. Hence, if we set the number of samples $R$ as above, then we have, with probability at least $1-(M+ M^2)\delta$:
    \begin{align}
        ||\delta F^R|| \leq \epsilon, \, ||\delta \mathcal{G}^R||\leq M\epsilon, 
    \end{align}
    where $||.||$ is the max-norm or the matrix norm induced by the max-norm depending on the argument being a vector or a matrix.\\
     Then, the R-sample solution, after ignoring the second and higher order terms of the small variables $\delta F^R$ and $\delta \mathcal{G}^R$, can be written as:
    \begin{align*}
        \hat{C}_M^R &= (\mathcal{G} + \delta \mathcal{G}^R)^{-1}(F + \delta F^R),\\
        &\approx \mathcal{G}^{-1} F + \mathcal{G}^{-1}\delta F^R - \mathcal{G}^{-1}\delta \mathcal{G}^R\mathcal{G}^{-1} F.
    \end{align*}
    Noting that $\hat{C}_M = \mathcal{G}^{-1} F$, it follows that:
    \begin{align*}
        \hat{C}_M^R - \hat{C}_M =  \mathcal{G}^{-1}\delta F^R - \mathcal{G}^{-1}\delta \mathcal{G}^R\mathcal{G}^{-1} F,
    \end{align*}
    and hence,
    \begin{align*}
        ||\hat{C}_M^R - \hat{C}_M|| &\leq ||\mathcal{G}^{-1}||||\delta F^R|| + ||\mathcal{G}^{-1}||^2||F|| ||\mathcal{\delta G}^R||,\\
        &\leq K_G\epsilon + K_G^2K_FM \epsilon \equiv \mathcal{C} \epsilon,
    \end{align*}
where $||\mathcal{G}^{-1}|| = K_G, \, ||F|| = K_F,\, \mathcal{C} = K_G + K_G^2 K_F M$.\\
However, from Proposition 1, we know that $||C_M^* - C_M|| \leq \kappa(\mathcal{G})||C_M^*|| \varepsilon$. Therefore, it follows that the error between the R-sample solution and the true solution satisfy:
\begin{align*}
    ||\hat{C}_M^R - C_M^*|| \leq \mathcal{C}\epsilon + \kappa(\mathcal{G})||C_M^*||\varepsilon.
\end{align*}
The development above can be summarized as the following result.\\
\begin{proposition}
    Let $R$ be chosen as in Eq. \ref{sample_C}. Then, with probability at least $1-(M+M^2)\delta$, the R-sample estimate $\hat{C}_M^R$ satisfies: \begin{align*}
    ||\hat{C}_M^R - C_M^*|| \leq \mathcal{C}\epsilon + \kappa(\mathcal{G})||C_M^*||\varepsilon,
    \end{align*}
    where  $||\mathcal{G}^{-1}|| = K_G, \, ||F|| = K_F,\, \mathcal{C} = K_G + K_G^2 K_F M$, $||.||$ is the max-norm when applied to vectors, and the max-induced matrix norm when applied to matrices. The error parameters $\varepsilon$ and $\epsilon$ were respectively defined in Lemma \ref{e-lemma} and Proposition \ref{prop2}.\\
    \end{proposition}
\begin{remark}
 The above analysis suggests that given sufficient samples, the data-based solution can be a good approximation of the true solution given sufficiently small $\epsilon \text{and} \,\varepsilon$. However, note that the sample complexity for the error to be below level $\epsilon$ with high confidence $(1-\delta)$ is critically dependent on the variances $\nu_i^2$ and $\nu_{ij}^2$. Furthermore, the error in the solution $||C_M^R - C_M^*||$ is dependent on the conditioning $\kappa(\mathcal{G})$ of the underlying Gram matrix. 
\end{remark}

\subsection{Nonlinear Case}
Given the nonlinear approximation architecture, one has to solve the least-squares problem 
\begin{equation*}
    \theta^* = \arg\min_\theta \frac{1}{2}||f-h(x,\theta)||^2,
\end{equation*}

where recall that $\theta = [\theta_1, \theta_2, \cdots \theta_n]'$ is an n-dimensional parameter.\\
This may be solved using a suitable optimization algorithm such as gradient descent. Let 
\begin{equation}
    J(\theta)=\frac{1}{2}||f-h(x,\theta)||^2= \frac{1}{2}\int (f-h(x,\theta))^2 \pi(x)dx.
\end{equation} This implies that 
\begin{equation*}
    \frac{\partial J}{\partial \theta} = -\int (f-h(x,\theta))h_\theta (x,\theta)\pi(x)dx,
\end{equation*}
where $h_\theta(x,\theta)=\frac{\partial h(x,\theta)}{\partial \theta}=[\frac{\partial h}{\partial \theta_1}, \cdots, \frac{\partial h}{\partial \theta_n}]'. $ 

Then, we can set up the gradient descent recursion 
\begin{equation*}
    \theta_{t+1} = \theta_t - \gamma_t \frac{\partial J}{\partial \theta}\Bigg|_{\theta_t},
\end{equation*}
where $\theta_t \rightarrow\theta^*$ s.t. $\frac{\partial J}{\partial \theta}|_{\theta^*}=0$, given $\gamma_t \rightarrow 0$ and $\sum_t \gamma_t = \infty$ (see Ch. 3, Prop. 3.5 in \cite{NDP}). Note, however, that, unlike the linear case, there is no guarantee that the algorithm converges to the global minimum.

Next, we consider the data-driven estimation problem for the nonlinear parametrization. Given the approximation architecture $h(x,\theta)$, parameterized non-linearly by $\theta$, 
let us define our data-based cost function as 
\begin{equation*}
    J^R(\theta)=\frac{1}{2R}\sum_{j=1}^R (f(x^j)-h(x^j,\theta))^2.
\end{equation*} Thus, we have the data-based cost-gradient 
\begin{equation*}
    \frac{\partial J^R(\theta)}{\partial \theta}=- \frac{1}{R}\sum_{j=1}^R(f(x^j)-h(x^j,\theta))\frac{\partial h(x^j,\theta)}{\partial \theta}.
\end{equation*}
From the Law of Large Numbers, we have $\frac{\partial J^R}{\partial \theta}\rightarrow \frac{\partial J}{\partial \theta}$ as $R\rightarrow \infty$. Thus, as the number of samples becomes large, the gradient estimate becomes arbitrarily accurate and the gradient descent algorithm is guaranteed to converge to a stationary point of the cost function. In particular, the following result holds. Let:
\begin{align*}
    \Delta(x) = f(x) - h(x,\theta), \; H^i(x) = \frac{\partial{h}}{\partial \theta_i},\\
    \mathrm{Var}[\Delta(x) H^i(x)] = \mu_i^2.
\end{align*}
\begin{proposition}
Given the step size parameters $\gamma_t$ satisfy $\sum_t \gamma_t = \infty, \; \sum_t \gamma_t^2 < \infty$ and the variances are uniformly bounded $\mu_i^2 < K< \infty$ for all $\theta$, the stochastic gradient descent sequence: $\theta_{t+1} = \theta_t - \gamma_t \frac{\partial J^R}{\partial \theta}|_{\theta_t}$ converges to $\bar{\theta}$ such that $\frac{\partial J}{\partial \theta}|_{\bar{\theta}} = 0$, for any finite $R$.
\end{proposition}
\begin{proof}
    Let $\mathcal{F}^t$ denote the history of the algorithm till time $t$. Write $\frac{\partial J^R}{\partial \theta}|_{\theta_t} = \frac{\partial J}{\partial \theta}|_{\theta_t} + w_t$. Then, it follows due to the unbiased nature of the estimate that $\mathrm{E}[w_t/ \mathcal{F}^t] = 0$.  and $\mathrm{E}[||w_t||^2/ \mathcal{F}_t] \leq \frac{nK}{R} < \infty$, uniformly for all $t$. Therefore, the sequence $\theta_t$ is guaranteed to converge to a stationary point of the function $J(.)$ (Prop. 4.1, \cite{NDP}).
\end{proof}

Further, we can give sample complexity bounds for the error in the gradient estimate analogous to the linear case. Let:
\begin{align*}
    \nabla_iJ = \frac{\partial J}{\partial \theta_i},\; \nabla_i^R J = \frac{\partial J^R}{\partial \theta_i}.
\end{align*}
Then, due to the Central limit theorem, given $R$ is is large enough: $\nabla_i^RJ \sim \mathcal{N}(\nabla_iJ, \frac{\mu_i^2}{R})$. \\
Thus, for $P(|\nabla_i^RJ - \nabla_i J|> \epsilon|) < \delta$, the number of samples required is:
\begin{align}\label{sample_grad}
    R_i = \mu_i^2 \frac{p^2(\delta)}{\epsilon^2}.
\end{align}
Hence, given $R = \max_i R_i$, with probability at least $(1-n\delta)$, $|\Delta_i^R J - \Delta_i J| < \epsilon$, for all $i=1,2,\cdots n$.\\

\begin{remark}
    There exists work that establishes the universal approximation properties of neural networks, both in the bounded depth and arbitrary width case \cite{NN_UA1,NN_UA2}, and more recently, in the bounded width and arbitrary depth case \cite{NN_UA3,NN_UA5,NN_UA4}. However, these are existence results and the actual training of these networks from data requires approaches such as backpropagation for which there are relatively few guarantees beyond convergence to a stationary point using stochastic gradient descent as outlined above.  The fundamental issue is that albeit there may be a unique parameter value that minimizes the cost, gradient descent cannot be expected to converge to that value. Further, as the size of the network is increased (either depth/ width), is there a unique or are there multiple parameter values that converge to the function being approximated, and if so, is there any guarantee that gradient descent via backpropagation can converge to any of these parameter values? These are open questions in the literature and beyond the scope of this paper, and thus, we do not provide a data-based approximation result such as Proposition 3 in the linear case. Nonetheless, in what follows, we will show that ill-conditioning and variance affect the training process adversely if approximating a function over a large domain/ high order, and consequently,  there is no guarantee regarding the convergence of the parameter estimates beyond the fact that they converge to a stationary point of the cost. Further, even with the same dataset, owing to the ill-conditioning and variance, gradient descent cannot be expected to converge to the same stationary point over multiple different runs.
\end{remark}
\section{Variance and Ill-Conditioning in the Model Estimation Problem}\label{sec_var_illcon}
In this Section, we shall show that the data-based solution suffers from very high variance as well as severe ill-conditioning if approximating over a large domain. 

In order to make an approximation of the function $f(.)$ over a large domain, one has to sample with a large enough exploration parameter $\sigma$ as well as expand to a high enough order $M$. This, in turn, leads to two issues:

\textbf{Variance.} $\mathrm{Var}(f\psi^i) = \nu_i^2$ and $\mathrm{Var}(\psi_i\psi_j) = \nu_{ij}^2$ grows combinatorially-exponentially in the order $M$ and the size of the sampled domain $\sigma$.\\
\textbf{Ill-Conditioning.} The elements of the Gram matrix $\mathcal{G}_{ij}$ are exponentially-combinatorially different, leading to severe ill-conditioning.

As mentioned previously, the data based solution may be written as:
\begin{equation}
    \hat{C}^R_M=(\underline{\Psi}^{R'}\underline{\Psi}^{R})^{-1}\underline{\Psi}^{R'}F^R = (\mathcal{G}+\delta G^R)^{-1}(F+\delta F^R), \label{ill-cond}
\end{equation} 
where $\mathcal{G}$ and $F$ are the true Gram matrix and projections of $f$ on the basis vectors respectively. If the approximation order $M$ is high, the variance in the errors $\delta G^R$ and $\delta F^R$ is very large, and $\mathcal{G}$ is ill-conditioned, and thus, the data-based solution $\hat{C}^R_M$ is bound to have a large variance, and also is sensitive to small variations owing to the ill-conditioning. In the following development, we will characterize the variance in the estimates for some common choices of sampling distributions and generic basis functions.\\

\begin{remark}
    We note that the ill-conditioning is a feature of the problem formulation and is due to the high order approximation over a large domain. In addition, the variance stems from the data-based solution to the problem. 
\end{remark}

\subsection{Linear Case}
We have seen that the optimal estimates $\hat{C}_M$ and $C_M^*$ are $\varepsilon$-close (Proposition 2).  
Furthermore, as $R \to \infty$, we have shown that the data-based estimate $\hat{C}_M^R \to \hat{C}_M$.

Next, we need to characterize the variance of the data-driven solution:
\begin{align}
    \mathcal{G}^R \, \hat{C}_M^R &= (\Psi^R)' F^R,
\end{align}
where
$
    (\Psi^R)' \Psi^R = \mathcal{G}^R.
$

Noting that 
\begin{align}
   [ (\Psi^R)' f^R]_i &\to \mathbb{E}[f \, \psi^i], \nonumber
\end{align}
we need to find the mean $f_i$ and the variance $\nu_i^2$ of the variable $f (x)\, \psi^i(x)$, where $x \sim \pi(.)$, i.e., $x$ is sampled from the distribution $\pi(.)$ .
Further,
\begin{align}
    (\Psi^R)' \Psi^R &\to \mathcal{G} = [ \gamma_{ij}(\sigma) ] = \left[ \mathbb{E}[\psi_i \, \psi_j] \right], \nonumber
\end{align}
so we need to find the mean $\gamma_{ij}$ and the variance $\nu_{ij}^2$ of the random variable $\psi_i(x)  \, \psi_j(x)$, $x\sim \pi(.)$.

Define:
\begin{align}
    \mathbb{E}\left[ x^k \right] &= B_k \, \sigma^k. \\
    \text{Note that } &\text{if the variable $x$ is Gaussian, i.e., } x \sim \mathcal{N}(0, \sigma^2): \nonumber \\
    B_k &= 
    \begin{cases}
        \left( \frac{k}{2} - 1 \right)!!, & \text{if $k$ even}, \\
        0, & \text{if $k$ odd}.
    \end{cases}
\end{align}

\begin{align}
\text{Uniform: } U(-\sigma, \sigma), \quad 
B_k =
\begin{cases}
    \frac{1}{k+1}, & \text{if $k$ odd}, \\
    0, & \text{if $k$ even}.
\end{cases}
\end{align}

\begin{align*}
\text{Also, let:} \quad
f(x) &= \sum_k \frac{\phi_k}{k!} x^k, \\
\psi^i(x) &= \sum_m \frac{\beta_m^i}{m!} x^m, 
\end{align*}
the Taylor expansions of the functions $f$ and $\psi^i$ respectively.
\begin{lemma} \quad 
The mean $f_i(\sigma)$ and the variance $\nu_i^2(\sigma)$ of the r.v. $f \psi^i$ is given by: 
\begin{align*}
\mathbb{E}[f\psi^i]&=f_i(\sigma) = \sum_{m,n} \frac{\beta^i_m \phi_n}{m! \, n!} \, B_{n+m} \, \sigma^{n+m}, \\[0.5em]
\mathrm{Var}(f \psi^i) &= \nu_i^2(\sigma) \\
&= \sum_{m,n,p,k} \frac{\phi_n \phi_k \beta^i_m \beta^i_p}{m! \, n! \, p! \, k!} \\
&\times \left( B_{n+m+p+k} - B_{n+m} B_{p+k} \right) \sigma^{n+m+p+k}.
\end{align*}
\end{lemma}

\begin{proof} 
\begin{align*}
f\psi^i = \sum_{m,n} \frac{\phi_n \beta^i_m}{m! \, n!} x^{n+m}. 
\end{align*}
Thus,
\begin{align*}
f_i(\sigma) = \mathbb{E}[f \psi^i] 
&= \sum_{m,n} \frac{\phi_n \beta^i_m}{m! \, n!} \mathbb{E}[x^{n+m}]\\
&= \sum_{m,n} \frac{\phi_n \beta^i_m}{m! \, n!} B_{n+m} \sigma^{n+m}.
\end{align*}
Similarly,
\begin{align*}
\nu_i^2 &= \mathrm{Var}(f \psi^i) 
= \mathbb{E} \big[ (f \psi^i - f^i)^2 \big] \\
&= \sum_{m,n,p,k} \frac{\phi_n \beta^i_m \phi_k \beta^i_p}{m! \, n! \, p! \, k!} 
\left[ \mathbb{E}[x^{n+m+p+k}] - \mathbb{E}[x^{n+m}] \mathbb{E}[x^{p+k}] \right] \\
&= \sum_{m,n,p,k} \frac{\phi_n \phi_k \beta^i_m \beta^i_p}{m! \, n! \, p! \, k!} 
\left( B_{n+m+p+k} - B_{n+m} B_{p+k} \right) \sigma^{n+m+p+k}.
\end{align*}
\end{proof}

By a similar development, we have the following result for the Gram matrix element  $(\psi_i \psi_j)$.\\
\begin{corollary} The mean $\gamma_{ij}(\sigma)$ and the variance $\nu_{ij}^2(\sigma)$ of the r.v. $\psi_i \psi_j$ are given by:
\begin{align*}
\gamma_{ij}(\sigma) &= \sum_{m,n} \frac{\beta^i_m \beta^j_n}{m! \, n!} \, B_{n+m} \, \sigma^{n+m}, \\[0.5em]
\nu_{ij}^2(\sigma)  &= \sum_{m,n,p,k} \frac{\beta^i_m \beta^j_n \beta^i_p \beta^j_k}{m! \, n! \, p! \, k!} \\
&\times\left( B_{n+m+p+k} - B_{n+m} B_{p+k} \right) \sigma^{n+m+p+k}.
\end{align*}
\end{corollary}

\textbf{Variance:}  
Note that the variance of $f \psi^i$ is an infinite sum over $m,n,p,k$. In general, this sum is difficult to evaluate. Thus, in order to get an idea about at least the order of the variance in an easily comprehensible form, assume that the maximum of the summand when $m=n=p=k$ occurs at some $N_i$. 
Then, using Lemma 2 and Corollary 1, it follows that:
\begin{align} \label{varianceF}
\nu_i^2=\mathrm{Var}(f \psi^i) \gtrapprox  O \left( \frac{B_{4N_i} - B_{2N_i}^2}{(N_i!)^4} \sigma^{4N_i} \right).
\end{align}
Similarly, the variance of $\psi_i\psi_j$ is an infinite sum and suppose that the maximum of the summand when $m=n=p=k$, occurs for some $L_{ij}$, then:
\begin{align} \label{varianceG}
\nu_{ij}^2=\mathrm{Var}(\psi^i \psi^j) \gtrapprox  O \left( \frac{B_{4L_{ij}} - B_{2L_{ij}}^2}{(L_{ij}!)^4} \sigma^{4L_{ij}} \right).
\end{align}

\textbf{Sample Complexity Revisited:}
    Using Proposition 3, it can be seen that the number of samples required such that $|f_i^R -f_i| > \epsilon$, where recall that $f_i^R$ is the R-sample average estimate of the mean $f_i$, with probability at most $\delta$ is given by:
    \begin{align}\label{Fsample}
        R_i \gtrapprox O\left(\frac{\nu_i^2}{\epsilon^2}p^2(\delta)\right) = O\left(\frac{B_{4N_i} - B_{2N_i}^2}{(N_i!)^4}\sigma^{4N_i} \frac{p^2(\delta)}{\epsilon^2} \right),
    \end{align}
    where $p(\delta) = \Phi^{-1}(1-\delta)$. For instance, for $\delta = 0.001$, $p(\delta) \approx 3$.
    Similarly, if we want each $|\gamma_{ij}^R - \gamma_{ij}| > \epsilon$ with probability at most $\delta$, where $\gamma_{ij}^R$ is the R-sample average estimate of the Gram matrix element $\gamma_{ij}$, then the sample complexity is given by:
    \begin{align}\label{Gsample}
        S_{ij} \gtrapprox O\left(\frac{B_{4L_{ij}} - B_{2L_{ij}}^2}{(L_{ij}!)^4}\sigma^{4L{ij}} \frac{p^2(\delta)}{\epsilon^2} \right).
    \end{align}
    Thus, the sample complexity required is:
    \begin{align}\label{sample_Capp}
        \mathcal{R} = \max \left( \max_i R_i, \max_{i,j} S_{ij} \right).
    \end{align}

    \begin{remark}\label{exp_variance}
    The general summand term in the infinite sum for $\nu_i^2$ (see Lemma 2), say $\nu_i^{2,m,n,p,k} \rightarrow 0$, as $m,n,p,k \rightarrow \infty$ (as well as the terms for the variance $\nu_{ij}^2$). This is necessary for the infinite sum to be convergent and can be shown to be the case using Stirling's formula. However, this does not preclude the variance from becoming very large for some intermediate $N_i$ (respectively $L_{ij}$ for $\nu_{ij}^2$). In particular, the variance is determined by the interplay between the exploration parameter $\sigma$ and the ``implicit" order of the approximation $N_i$/ $L_{ij}$, which is in turn decided by the function that we are trying to approximate. Also, note that we are simply showing the peak variance term in the expressions above and the actual variance will be significantly larger since it is the infinite sum of the individual terms, and thus, the sample complexity expression above is a  rather crude lower bound. Nonetheless, the development above clearly shows that the sample complexity of the solution grows combinatorially-exponentially in the size $\sigma$ of the approximation domain and the order of the approximation required $(N_i, L_{ij})$. \\
\end{remark}
\textbf{Ill-conditioning:}  
Consider the term $\gamma_{ij}$, and suppose, as above, that the  maximum of the summand is attained when $m=n$ at $K_{ij}$.  
Then,
\begin{align} 
\gamma_{ij} &\gtrapprox O\!\left( \frac{B_{2K_{ij}} \, \sigma^{2K_{ij}}}{(K_{ij}!)^2} \right).
\end{align}
However, it is difficult to understand the ill-conditioning issue without an additional assumption. Suppose now that for any $\gamma_{ij}$, the summation is from $m=0:i, n= 0:j$. This implies that the basis functions form progressively higher order approximations, which is typically the case and is  borne out by the empirical evidence in Section \ref{sec4}. Then, note that $\gamma_{00} \approx O(1)$, and assuming that the maximum value of $\gamma_{ij}$ over all possible $(i,j)$, when $i=j$, is attained at some $\eta$, and defining $K = K_{\eta\eta}$, an approximate lower bound on the conditioning of the Gram matrix is given as:
\begin{align}\label{conditioning}
    \kappa(\mathcal{G}) \gtrapprox O\!\left( \frac{B_{2K}\sigma^{2K}}{(K!)^2} \right)
\end{align}

If we want to approximate the function $f(.)$ over a suitably large domain, then the exploration parameter $\sigma$ needs to be large as should be the order of approximations $K$ as defined above. Thus, it follows that the Gram matrix $\mathcal{G}$ is highly ill-conditioned, and combined with the variance in the estimates $f_i^R$, the solutions are subject to high variance and consequently unreliable (see Eq. \eqref{ill-cond}).\\

\begin{remark}
The error parameter $\epsilon$ and the confidence level $\delta$ play an identical role here as they do in Probably Approximately Correct (PAC) learning theory \cite{PAC-1,PAC-2}. However, the sample complexity is seen to be combinatorial-exponential in the size of the hypothesis class, in this case, the size of the approximation domain $\sigma$ and the effective order of the approximation $(N_i, L_{ij})$ (see equations \ref{Fsample} and \ref{Gsample} above). Thus, this is not an efficiently learnable problem in the PAC sense since in such a case the dependence has to be polynomial.
\end{remark}


\subsection{Nonlinear Case}\label{nonlinear_theory}
In the following, we show the variance inherent in the gradient estimation and the ill-conditioning of the nonlinear case, the development is almost identical to the linear case.\\
\textbf{Variance:}
The gradient term  in the nonlinear case is given by:
\[
\frac{\partial J}{\partial \theta}
= \mathbb{E} \left[ \Delta(x) \, H(x, \theta) \right],
\]
where
\[
H(x, \theta) =
\begin{bmatrix}
H^{1}(x,\theta) \\
\vdots\\
H^{n}(x,\theta)
\end{bmatrix}
=
\begin{bmatrix}
\frac{\partial h}{\partial \theta_1} \\
\vdots\\
\frac{\partial h}{\partial \theta_p}
\end{bmatrix},
\]
and
\[
\Delta(x) = f(x) - h(x, \theta).
\]

Let
\[
\Delta(x) = \sum_{k=0}^{\infty} \frac{\Delta_k}{k!} x^k,
\quad
H^{i}(x) = \sum_{k=0}^{\infty} \frac{h^{i}_k}{k!} x^k.
\]

Then, following the same development as in the linear case,  
the variance in the estimate of the gradient is given by:
\begin{align*}
&\mu_i^2 = \mathrm{Var}(\Delta H^{i}) =\\
\sum_{m,n,r,s}
\frac{\, h^{i}_m \, \Delta_n \, h^{i}_r \, \Delta_s}
{m! \, n! \, r! \, s!}&\left( B_{m+n+r+s} - B_{m+n} \, B_{r+s} \right) \sigma^{m+n+r+s}.
\end{align*}
In the expression above, given $i$, if the maximum value, when $m=n=p=k$, is attained at $N_i$, the variance is of the order of:
\[
\mu_i^2
\gtrapprox
O\!\left(
\frac{B_{4N_i} - B_{2N_i}^2}{(N_i!)^4} \,
\sigma^{4N_i}
\right).
\]
Then, an approximate lower bound on the number of samples required for a given error-confidence level of $(\epsilon, \delta)$ in the gradient estimate with respect to the $i^{th}$ parameter, $\nabla_i^R J$ is given by:
\begin{align}\label{samplecomp}
    R_i \gtrapprox O\left(\frac{B_{4N_i} - B_{2N_i}^2}{(N_i!)^4} \sigma^{4N_i}\frac{p^2(\delta)}{\epsilon^2} \right).
\end{align}
Thus, the sample complexity required: $R = \max_i R_i$, for the error in all the components of the gradient estimate to be within $\epsilon$ with probability at least $1-n\delta$. Hence, it is clear that even in the nonlinear case, the sample complexity of the gradient estimate is exponential-combinatorial in the parameters $(\sigma, N_i)$. \\

\textbf{Ill-conditioning:} In order to show the effect of ill-conditioning in the nonlinear case, consider the following. Let us suppose that the optimal parameter $\theta^*$ is close enough to the current estimate $\Bar{\theta}$ such that a linear expansion about $\Bar{\theta}$, in terms of the perturbation $\delta \theta = \theta - \bar{\theta}$, suffices to approximate $h(x,\theta)$, \textit{i.e.},
\begin{align*}
     &h(x,\theta)\approx h(x,\Bar{\theta})+\underbrace{\begin{bmatrix}
         \frac{\partial h(x,\theta)}{\partial \theta^1} &\cdots &\frac{\partial h(x,\theta)}{\partial \theta^n}
     \end{bmatrix} \Bigg|_{\bar{\theta}} }_{H(x,\bar{\theta})\equiv H} \underbrace{\begin{bmatrix}
         \delta \theta^1 \\ \vdots \\ \delta \theta^n
     \end{bmatrix}}_{\delta {\theta}}. \\
     &\mbox{Define:\;} \epsilon(\theta)= \underbrace{\begin{bmatrix}
         f(x^1)-h(x^1,\bar{\theta})\\ \vdots \\ f(x^R)-h(x^R,\bar{\theta})
     \end{bmatrix}}_{\bar{\Delta}^R} - \underbrace{ \begin{bmatrix}
         H(x^1,\bar{\theta})\\ \vdots \\ H(x^R,\bar{\theta})
     \end{bmatrix} }_{\bar{H}^R} \delta {\theta} 
 \end{align*}

We can rewrite the original cost function in terms of the perturbation variable as: $J(\theta)=\epsilon(\theta)^T \epsilon(\theta)=\epsilon(\delta \theta)^T \epsilon (\delta \theta)\equiv
J(\delta \theta)$. Therefore,
 \begin{equation*}
     J(\delta \theta)=(\Bar{H}^R \delta \theta - \Bar{\Delta}^R)^T (\Bar{H}^R \delta \theta - \Bar{\Delta}^R).
\end{equation*}

Then, the gradient descent algorithm in the vicinity of $\Bar{\theta}$ can be written as:
\begin{align*}
  \delta \theta_t &= \delta \theta_{t-1} - \gamma \frac{\partial J(\delta \theta)}{\partial (\delta \theta)} \\
  &= \delta \theta_{t-1} - \gamma \bar{H}^{R'} (\bar{H}^R \delta\theta-\bar{\Delta}^R),
\end{align*}
which, for small $\gamma$, can be approximated as $\delta \Dot{\theta}=-\bar{H}^{R'}(\bar{H}^R\delta\theta- \bar{\Delta}^R)$. At convergence, $\delta\Dot{\theta}=0$, which implies
\begin{equation}
\label{eq: LSEquiv}
   \delta \theta = (\bar{H}^{R'}\bar{H}^R)^{-1} (\bar{H}^{R'}\bar{\Delta}^R).
\end{equation}
Note that the above equation has the same least square solution form as the linear case and the empirical Gram matrix in this case is $\mathcal{G}^R = \bar{H}^{R'}\bar{H}^R$, i.e., the matrix $H^R$ containing the derivatives of the nonlinear approximator with respect to the parameters plays the role of the matrix $\Psi^R$ in the linear case, and will have the same ill-conditioning as the linear case.

\subsection{Multidimensional Case}
In the following, we show how to extend the variance analysis to a multidimensional function. We examine the case of a scalar valued multidimensional function, the extension to the vector valued case is relatively straightforward. Let $\mathbf{x}= (x_1, \cdots x_d) $. Further, let
\begin{align}
f(\mathbf{x}) &= \sum_{\substack{i_1+i_2+\cdots+i_d=n \\ n}}
   \frac{\phi_{i_1,\dots,i_d}}{i_1! i_2! \cdots i_d!}
   x_1^{i_1} x_2^{i_2} \cdots x_d^{i_d}, \\[1em]
\psi^i(\mathbf{x}) &= \sum_{\substack{i_1+i_2+\cdots+i_d=n \\ n}}
   \frac{\beta^i_{i_1,\dots,i_d}}{i_1! i_2! \cdots i_d!}
   x_1^{i_1} x_2^{i_2} \cdots x_d^{i_d},
\end{align}
denote the multidimensional Taylor series expansions of the function and the basis functions respectively. Then,
\begin{align}
f \psi^i 
&= \sum_{\substack{i_1+i_2+\cdots+i_d = m \\ k_1+\cdots+k_d = n \\
m,n}}
   \frac{\phi_{i_1,\dots,i_d} \, \beta^i_{k_1,\dots,k_d}}
        {i_1! \cdots i_d! \; k_1! \cdots k_d!}
   x_1^{i_1+k_1} \cdots x_d^{i_d+k_d}.
\end{align}
 To simplify the notation, let $I=\{i_1,\cdots i_d\}$. Further, we will write: $|I| = m$, to denote the sum $\sum_{k=1}^d i_k = m$ and $I!$ to denote the product $\prod_{k=1}^d i_k$ for all possible $I$.\\
 
Next, we assume that $x_1,\dots,x_d$ are sampled independently and each variable is sampled with the same exploration parameter $\sigma$. Then, we obtain the following.
\begin{align}
\mathbb{E}\!\left[ f \psi^i \right] &=f_i(\sigma) \\
&= \sum_{\substack{|I|=m \\ |K| = n \\ m,n}}
   \frac{\phi_I \,
   \beta^i_K \,}{I!K!}
   B_{i_1+k_1}\sigma^{i_1+k_1} \cdots B_{i_d+k_d}\sigma^{i_d+k_d}.
\end{align}

Next, we can write:
\[
\sigma_i^2=\mathrm{Var}(f \psi^i) 
   = \mathbb{E}\!\left[
      \left(
      \sum_{\substack{|I|=m \\ |K|=n\\ m,n}}
      \tfrac{\phi_{I}\,\beta^i_{K}}
            {I! \; K!}
      x_1^{i_1+k_1}\cdots x_d^{i_d+k_d} - f^i
      \right)^{\!2}
   \right].
\]

Expanding,
\[
\begin{aligned}
\sigma_i^2
   &= \sum_{\substack{|I|=m,\, |K|=n \\ 
                      |L|=r,\, |P|=s \\ m,n,r,s}}
   \tfrac{\phi_{I}\,\phi_{L}\,
          \beta^i_{K}\,\beta^i_{P}}
        {(I!)(K!)(L!)(P!)} \\[0.5em]
   &\quad \times \Bigl[\prod_{r=1}^d B_{i_r+k_r+l_r+p_r} 
       - \prod_{r=1}^d B_{i_r+k_r}\,B_{l_r+p_r}\Bigr]\,
   \sigma^{m+n+r+s}.
\end{aligned}
\]
Suppose now that the maximum value when $m=n=r=s$ is attained at $N_i$, then each term in the sum above is of
$O\!\left( \dfrac{B_{4N_i} - B_{2N_i}^2}{(N_i!)^4} \, \sigma^{4N_i} \right)$. In particular, the value $\dfrac{B_{4N_i} - B_{2N_i}^2}{(N_i!)^4}$ in the $O(.)$ calculation is attained when $i_r=k_r=l_r=p_r=N_i$ for any $r=1,\cdots d$. However, evaluating the other terms, when more than one of the indices is non-zero is intractable. Instead, we substitute the $O(.)$ calculation above and note that this is simply to get a rough idea of the numbers involved.

Since there are $d^{N_i}$ ways in which each of $|I| = N_i$, $|K| = N_i$, $|L| = N_i$, and $|P| = N_i$, there are $d^{4N_i}$ such terms in the sum above.
Hence, it follows that
\[
\sigma_i^2 \gtrapprox O\!\left( \frac{B_{4N_i} - B_{2N_i}^2}{(N_i!)^4} \, (\sigma d)^{4N_i} \right).
\]
An analogous argument shows that $\sigma_{ij}^2 = \mathrm{Var}(\psi^i\psi^j)$ is also of the same order:
\begin{align}
    \sigma_{ij}^2 \gtrapprox O\!\left( \frac{B_{4L_{ij}} - B_{2L_{ij}}^2}{(L_{ij}!)^4} \, (\sigma d)^{4L_{ij}} \right),
\end{align}
where $L_{ij}$ is defined in an identical manner to the scalar case. The sample complexity is then given by equations \ref{Fsample}, \ref{Gsample} and \ref{sample_Capp}.

Finally, the conditioning of the problem is given by: 
\[
\kappa(\mathcal{G}) \gtrapprox O\!\left( \frac{B_{2K}}{(K!)^2} \, (\sigma d)^{2K} \right),
\]
with an analogous definition of $K$ as in the scalar case.
 We see that the dimension affects the variance and conditioning of the problem by effectively amplifying the exploration parameter $d$ times. \\
 \begin{remark}
 A note is due regarding the independent sampling in the different dimensions. It is always possible to linearly transform any set of correlated random variables such that they are uncorrelated and have the same variance, and thus, the variance expression above will still hold albeit the variables are uncorrelated, in general, and not necessarily independent. The resulting variables become independent if the sampling is Gaussian.\\
 \end{remark}
 
 \begin{remark}
     The extension of the multidimensional variance and ill-conditioning analysis to the nonlinear case can be done in a relatively straightforward fashion using the analysis in the scalar case and is left out for brevity.
 \end{remark}
\section{Recursive Linear Time-Varying Models and Nonlinear Optimal Control}\label{section5}
In this Section, we return to our original goal of solving a nonlinear optimal control problem. 
Thus far, we have seen that if we want a high-order/global approximation of a function $f(x)$, by sampling from a given distribution $\pi(\cdot)$, we are bound to suffer from high variance and ill-conditioning if approximating over a large domain. Consequently, a global data-driven surrogate dynamics model is bound to be inaccurate, and solving the optimal control problem utilizing such a one-shot identify-then-control surrogate model necessarily leads to an inaccurate solution. 
In the following, we show that if we ``limit" ourselves to estimating local LTV models, then the variance and ill-conditioning inherent to global model estimation cease to be a factor. Nonetheless, in spite of its local nature, when applied recursively, this approach can be used to accurately solve the motivating nonlinear optimal control problem.

\subsection{Iterative Linear Quadratic Regulator (ILQR)}\label{ILQR}
Let us recall the optimal control problem (Sec.~\ref{probform}). 
The ILQR Algorithm is a sequential quadratic programming technique for optimal control that \cite{ilqg2,li2004iterative,SQP_OC,SQP_OC-ACC26}, iteratively solves the nonlinear optimal control problem  using the following steps:

\subsubsection*{Step 1: Forward  Pass}
Given a nominal control sequence $\{u_t\}_{t=0}^{T-1}$ and initial state vector $x_0$, the state is propagated in time in accordance with the dynamics.
$
    x_{t+1}=f(x_t, u_t)
$
.Thus, we get the nominal trajectory $(\bar{x}_t,\bar{u}_t)$ and the linearized system around the nominal:

$
    \delta x_{t+1} = A_t \delta x_t + B_t \delta u_t.
$

\subsubsection*{Step 2: Backward pass}
Using the linearized system, the ILQR algorithm computes a local optimal control by solving the following discrete time Riccati Equation, which is given by,
\begin{align}
    \delta u_t
    &=R^{-1}B_t^T(-v_{t+1}-V_{t+1}(A_t\delta x_t + B_t\delta u_t))-\bar{u}_t, \nonumber \\
    &=-(R+B^T_t V_{t+1}B_t)^{-1}(R\bar{u}_t+B_t^T v_{t+1}\nonumber \\
    &+B_t^T V_{t+1}A_t\delta x_t),
\end{align}
which can be written in the linear feedback form 
\begin{equation*}
    \delta u_t =-k_t-K_t\delta x_t,
\end{equation*}
where 
\begin{align*}
    k_t &=(R+B_t^T V_{t+1}B_t)^{-1}(R\bar{u}_t+B_t^T v_{t+1}), \text{and}\\
    K_t &=(R+B_t^T V_{t+1}B_t)^{-1}B^T_t V_{t+1}A_t,
\end{align*}
with
\begin{subequations}
\begin{align}
    v_t &= l_{t,x}+A_t^T v_{t+1}-A_t^T V_{t+1}B_t(R+B_t^T V_{t+1}B_t)^{-1}\nonumber \\
    &\cdot(B_t^T v_{t+1}+R\bar{u}_t) \label{eq:vt}\\
    V_t &= l_{t,xx}+A_t^T (V_{t+1}^{-1}+B_t R^{-1}B_t^T)^{-1}A_t\nonumber \\
    &=l_{t,xx}+A_t^T V_{t+1}A_t-A_t^T V_{t+1}B_t(R+B_t^T V_{t+1}B_t)^{-1} \label{eq:Vt} \nonumber \\
    &\cdot B_t^T V_{t+1}A_t.
\end{align}
\end{subequations}
with the terminal conditions $v_T(x_T)=\frac{\partial c_T}{\partial x}|_{x_T}=\frac{\partial l}{\partial x}|_{x_T}$ and $V_T(x_T)=\nabla ^2_{xx}c_T|_{x_T}$. To compute $v_t$ and $V_t$, we need their corresponding values in the future, \textit{i.e.} $v_{t+1}$ and $V_{t+1}$. From Eqs.~\ref{eq:vt} and \ref{eq:Vt} we see that, given the terminal conditions and the local LTV model parameters $(A_t, B_t)$, we can do a backward-in-time sweep to compute all values of $v_t$ and $V_t$, and thus, the corresponding optimal control for that trajectory.
\subsubsection*{Step 3: Update trajectory}
Now, given the gains from the backward pass, we can update the nominal control sequence as
\begin{align}
    \bar{u}_t^{k+1} &= \bar{u}_t^k + k_t + K_t(x_t^{k+1}-x_t^k) \text{,} \nonumber\\
    x^{k+1}_0 &= x^k_0. \label{ilqr_update}
\end{align}
If the convergence criteria is not met ($k_t \approx 0$), we iterate the process by repeating steps (1) to (3).

\subsection{Recursive LTV Model Estimation and Optimal Control}

The key to solving the optimal control problem via the  ILQR approach is the LTV system $(A_t,B_t)$ found by linearizing the system around the current estimate of the optimal trajectory. However, this LTV system can be estimated from data as follows.\\
Assume that $x_0$ is given to us, as is an initial guess for the optimal control (say $\{u_t^0\}$). Suppose that we perturb the system as $x_0+\delta x_0$, and $\bar{u}_t^0+\delta u_t$, where $\delta x_0 \sim p_0^x(\cdot)$ and $\delta u_t \sim p_t^u(\cdot)$, and the sampling distributions are small enough ($\sigma$) such that the system is arbitrarily well approximated as linear, \textit{i.e.}, $\delta x_{t+1}\approx A_{t}\delta x_{t}+B_{t}\delta u_{t}$. Note that this is always feasible owing to Proposition \ref{ME_acc}.
Then, one may perform roll-outs to get sample trajectories:\\
\begin{align}
\delta x_{t+1}^i &\approx A_{t}\delta x_{t}^i+B_{t}\delta u_{t}^i \nonumber\\
\implies \underbrace{\begin{bmatrix}
        \delta x_{t+1}^1 \\ \vdots \\ \delta x_{t+1}^R 
    \end{bmatrix}}_{\delta x_{t+1}^R} &= \underbrace{\begin{bmatrix}
        \delta x_{t}^1  &\delta u_t^1 \\ \vdots &\vdots \\ \delta x_{t}^R &\delta u_t^R
    \end{bmatrix}}_{\delta X_{t}^R} \begin{bmatrix}
        A_t \\ B_t
    \end{bmatrix} \nonumber \\
    \implies \begin{bmatrix}
        \Hat{A}_t \\ \Hat{B}_t
    \end{bmatrix} &= (\delta X_t^{R'}\delta X_t^R)^{-1} (\delta X_t^{R'}\delta x_{t+1}^R). 
    \label{eq:LTV-LS}
    \end{align}
The empirical Gram matrix $(\delta X_t^{R'}\delta X_t^R)$ is well-conditioned as $N=1$ and $\sigma$ is small (see \eqref{conditioning}), and the variance in $(\delta X_t^{R'}\delta x_{t+1}^R)$ is low for the same reason (see \eqref{varianceF} and \eqref{varianceG}). Thus, we can estimate the LTV system $(A_t,B_t)$ around the trajectory generated from $(x_0,\bar{u}_t^0)$ in an accurate fashion.
However, this only gives us a local description of the system around the initial trajectory $(x_t^0, \bar{u}_t^0)$. Thus, one may question the practicality of such a local approximation. \\
In fact, this local description $(A_t, B_t)$ around the trajectory is all that is required: the local LTV estimate can be used to compute an improved control sequence, which in turn defines a new trajectory around which a new LTV system can be identified, and so on until convergence. This is precisely the `local' model-based RL approach developed by the authors and their co-workers, called the Decoupled Data-based Control (D2C), which proceeds as follows \cite{D2C_ASME}:
\subsubsection*{Step 1: Forward Pass}
The nominal control sequence $\{\bar{u}_t^k\}$ is propagated through the black-box nonlinear dynamics starting from $x_0$ to compute the nominal trajectory $(\bar{x}_t^k, \bar{u}_t^k)$.

\subsubsection*{Step 2: Data-based LTV System Identification}
Around the nominal trajectory $(\bar{x}_t^k, \bar{u}_t^k)$, the system is perturbed and $R$ rollouts are performed to estimate the local LTV model $(\hat{A}_t^k, \hat{B}_t^k)$ via the least-squares procedure outlined above (Eq.~\ref{eq:LTV-LS}).

\subsubsection*{Step 3: Backward Pass}
Given the terminal conditions $v_T(x_T) = \frac{\partial c_T}{\partial x}|_{x_T}$ and $V_T(x_T) = \nabla^2_{xx}c_T|_{x_T}$, along with the estimated LTV parameters $(\hat{A}_t^k, \hat{B}_t^k)$, a backward-in-time sweep is performed via Eqs.~\eqref{eq:vt} and \eqref{eq:Vt} to compute $\{v_t, V_t\}$ at all time steps, and thereby the feedforward and feedback gains $\{k_t^k, K_t^k\}$.

\subsubsection*{Step 4: Update Trajectory}
Given the gains from the backward pass, the nominal control sequence is updated according to Eq.~\eqref{ilqr_update}.

If the convergence criterion ($k_t^k \approx 0$) is not met, Steps 1--4 are iteratively repeated with the refined control sequence (Fig.~\ref{fig:optcontrol}(a)-(b)). At convergence, the procedure yields an optimal nominal control sequence $\bar{u}_t^*$ together with a local LTV perturbation feedback gain $K_t^*$. This shows that local LTV models suffice for optimal control owing to their local nature, thereby alleviating the ill-conditioning and variance inherent in estimating ``global or higher-order'' models in a data-based fashion.\\
\begin{remark}
Convergence of the LTV-ILQR Approach. More details regarding this approach can be found in the reference \cite{D2C_ASME} including the convergence of the approach starting from any random initial guess and the effect of identification errors on the convergence. In particular, it is shown that the approach converges regardless of the initial guess and errors in the system identification, given the errors are small enough.\\
\end{remark}

\begin{remark}
Classical Adaptive Control \cite{kumar2015stochastic,ioannou2012robust}. Typically, classical adaptive control is tailored for Linear Time Invariant (LTI) problems. LTI models describe the behavior of a system around an equilibrium, in contrast to nonlinear problems considered here where the system is far from an equilibrium, which lead to LTV problems when linearizing about a trajectory as shown above (an equilibrium is a trajectory that does not change). The LTV-ILQR approach can be thought of as a generalization of LTI based classical adaptive control to ``general" nonlinear problems. An important distinction is that adaptive control problems can be solved on a single sample path/ rollout owing to their LTI nature whereas the LTV approach requires ensembles of sample paths/ rollouts in order to solve the nonlinear problem. Indeed, this is the fundamental nature of the nonlinear problem.
\end{remark}

\subsection{Necessity of Recursive Local Modeling for Optimal Control} 
We have shown above that the LTV estimation approach suffices to solve the optimal control problem. However, the LTV approach is not necessary. For instance, instead of ILQR, if we use the related algorithm called differential Dynamic Programming (DDP) to find the optimal control \cite{ddp}, we need the Hessian information around the trajectories in addition to the Jacobean/ LTV information. Thus, in this case, we need to estimate the first and second order perturbation terms simultaneously which is still feasible since the order of approximation, $N=2$, is low. It is thought that such second order information can lead to faster convergence to the optimum albeit our recent work shows that using this second order information can actually adversely affect convergence when far from the optimum \cite{SQP_OC,SQP_OC-ACC26}. In practice, we have seen that estimating the Hessian terms simultaneously with the LTV terms can be prone to errors (owing to the variance and ill-conditioning outlined in this paper) and estimating only the LTV terms is more accurate. Nonetheless, we cannot theoretically claim that the LTV approach is necessary since moderate higher order information can be incorporated and perhaps lead to a better convergence behavior when compared to the LTV approach. However, any such method will still need to \emph{recursively identify the local (not necessarily LTV)  model around the current trajectory as the optimization proceeds} since the one-shot global approach to model estimation is bound to suffer from the Curse of Variance and Ill-Conditioning (COVIC) and cannot be used for control design.


      


%% file: Simulation_new.tex
\section{Empirical Results}\label{sec4}
In this section, we present the empirical results to validate our theoretical claims. The algorithms are trained on the MUJOCO-based simulators for the Pendulum, Cartpole and Fish systems, and bench-marked to compare performance. Furthermore, we show the discrepancy between models when training on local vs global neighbourhoods, and empirically validate the theoretical development in the prior sections.



\subsection{Sampling Distributions for Learning Models}
Using the optimal control algorithm (ILQR) as described in Sec \ref{ILQR}, we generate the optimal trajectory and control for the modified pendulum swingup task: Starting from the downward position of the pendulum ($\theta_0=180^\circ$), the pendulum needs to swing-up and come to rest at $\theta_T = 60^\circ$ within 30 time-steps. We generate the training data by increasingly perturbing the optimal control actions about the optimal trajectory, and sampling 2000 trajectories (60000 data-points) from the neighborhoods obtained. The sampled trajectories are then split in a 9:1 ratio for the testing and training set respectively. 

\begin{figure}[!htpb]
      \subfloat[Swingup with $10\%$ perturbation]{\includegraphics[width=0.5\linewidth]{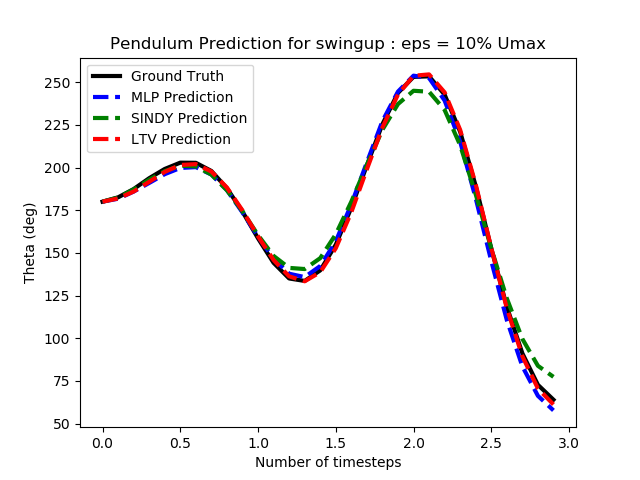}} 
      \subfloat[Phase-phase plot for $10\%$ perturbation]{\includegraphics[width=0.5\linewidth]{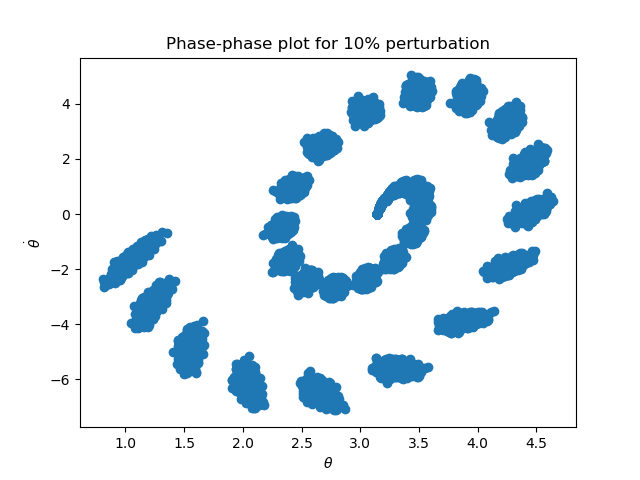}}
      
      \subfloat[Swingup with $30\%$ perturbation]{\includegraphics[width=0.5\linewidth]{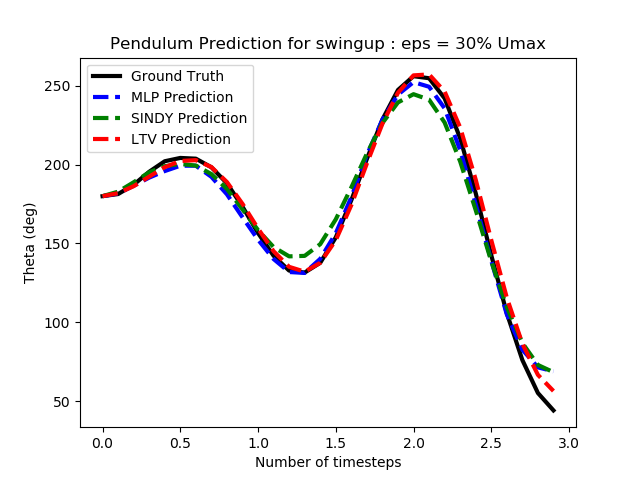}}
      \subfloat[Phase-phase plot for $30\%$ perturbation]{\includegraphics[width=0.5\linewidth]{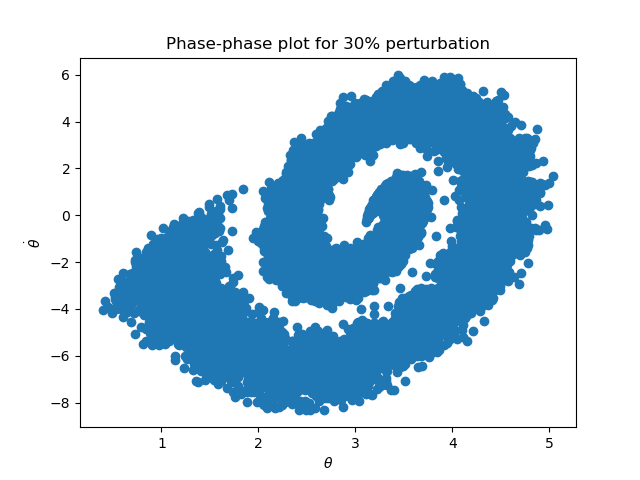}}
      
      \subfloat[Swingup with $60\%$ perturbation]{\includegraphics[width=0.5\linewidth]{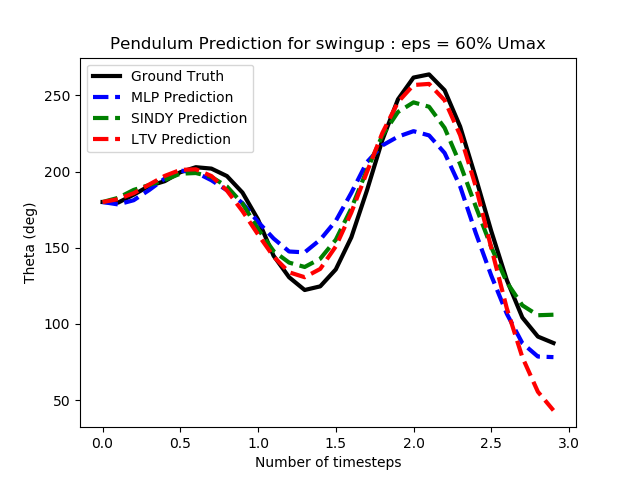}}
      \subfloat[Phase-phase plot for $60\%$ perturbation]{\includegraphics[width=0.5\linewidth]{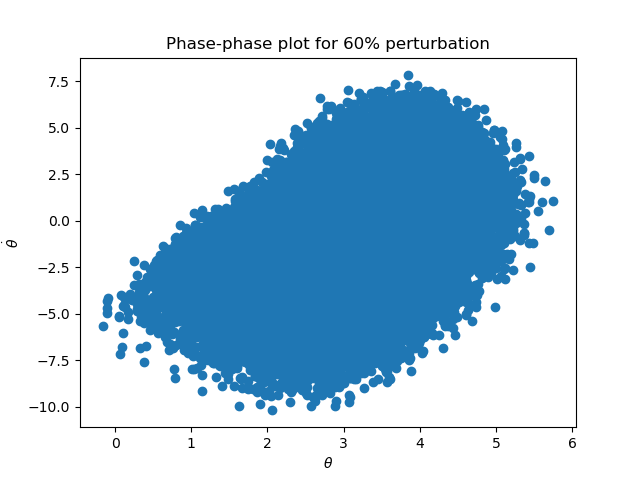}}
      

\caption{{Comparing the performance of the nonlinear SINDy and MLP models with the LTV model on swing-up pendulum trajectories.}}
\label{fig:LTV}
\end{figure}

As an analogue to the linear (in parameter) approximator, we make use of the established SINDy algorithm \cite{sindy}, directly adopting the open-source Python implementation PySINDy \cite{pysindy} for our comparison.
We train both the MLP model and SINDy on these data-points and generate the predictions to compare their performance. To benchmark these `global' nonlinear models against a `local' model, we compare them with an LTV model (Fig.~\ref{fig:LTV}). 

 Please note that the random perturbations of the optimal control implicitly define the sampling distribution in the state space (Fig.~\ref{fig:LTV}).

Similar to the pendulum, we generate the optimal trajectory and control for the modified Cartpole swingup task: Starting from the downward position of the pendulum ($\theta_0=180^\circ$) and the cart at rest, the cart needs to move such that the pendulum comes to rest at $\theta_T = 60^\circ$ within 30 time-steps. The training and testing data size and sampling methods remain the same. 

\begin{figure}[!htpb]
\begin{multicols}{2}
      \subfloat[$5\%$ Perturbation]{\includegraphics[width=\linewidth]{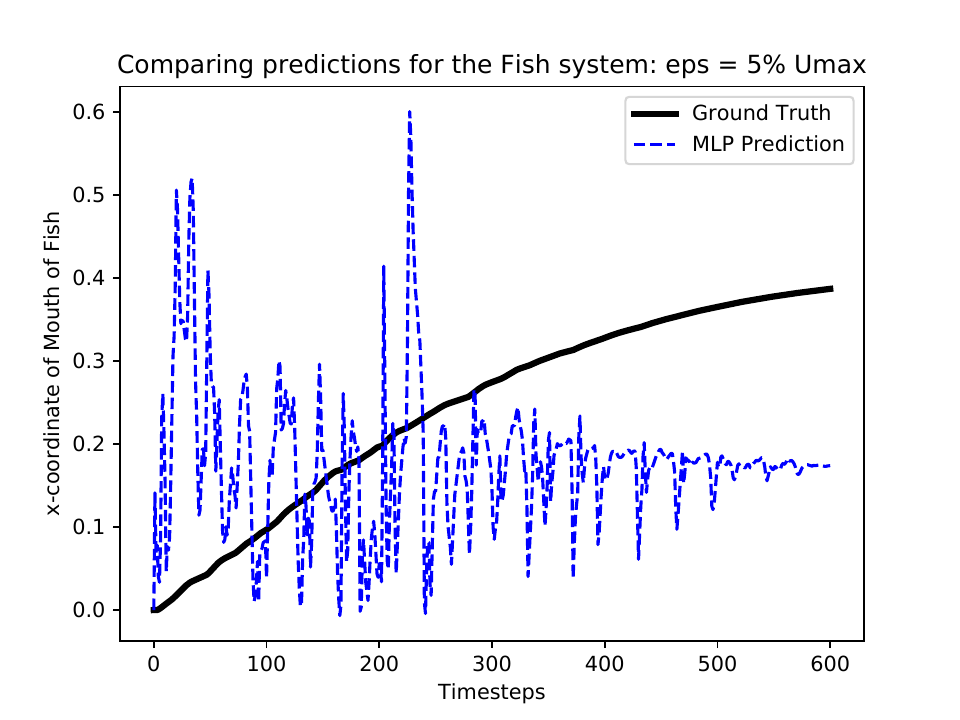}}    
      \subfloat[$30\%$ Perturbation]{\includegraphics[width=1\linewidth]{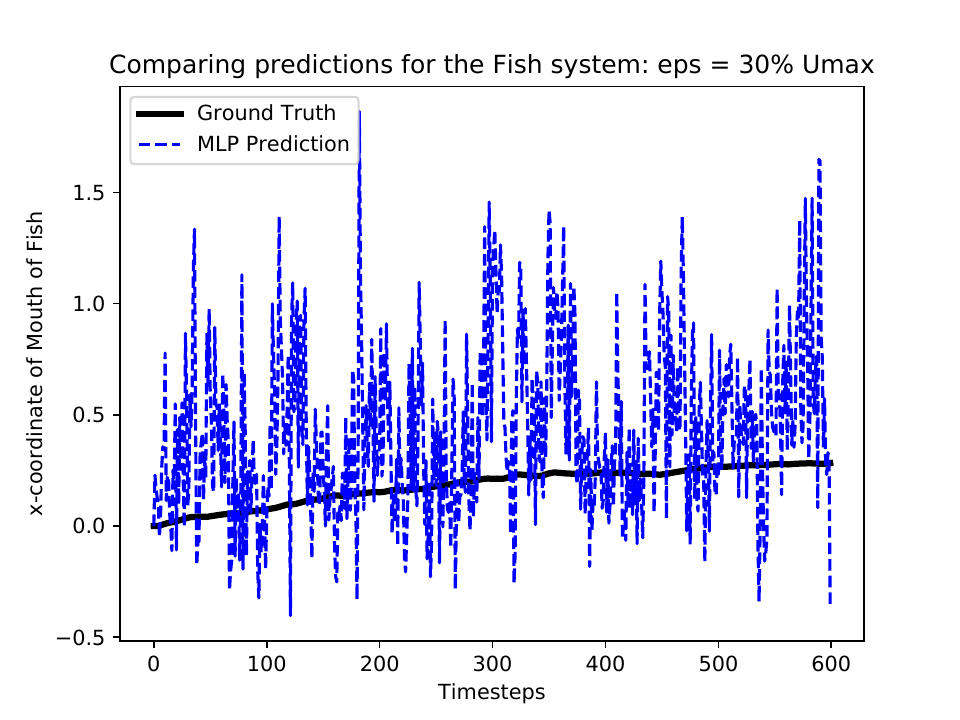}}
\end{multicols}
\begin{multicols}{2}
      \subfloat[$5\%$ Perturbation]{\includegraphics[width=\linewidth]{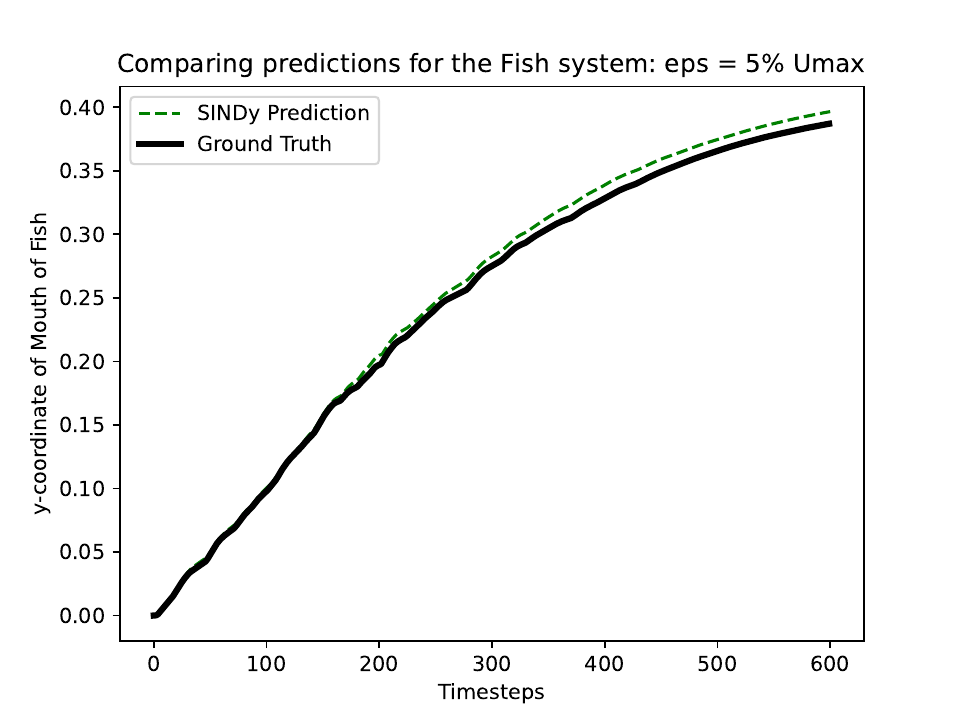}}    
      \subfloat[$30\%$ Perturbation]{\includegraphics[width=1\linewidth]{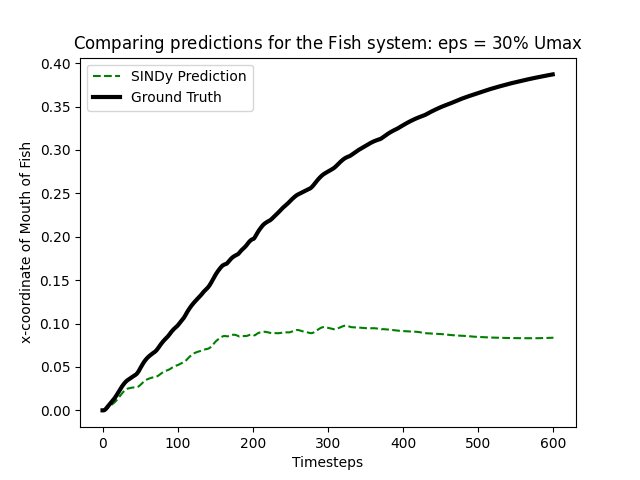}}
\end{multicols}
\begin{multicols}{2}
      \subfloat[$5\%$ Perturbation]{\includegraphics[width=\linewidth]{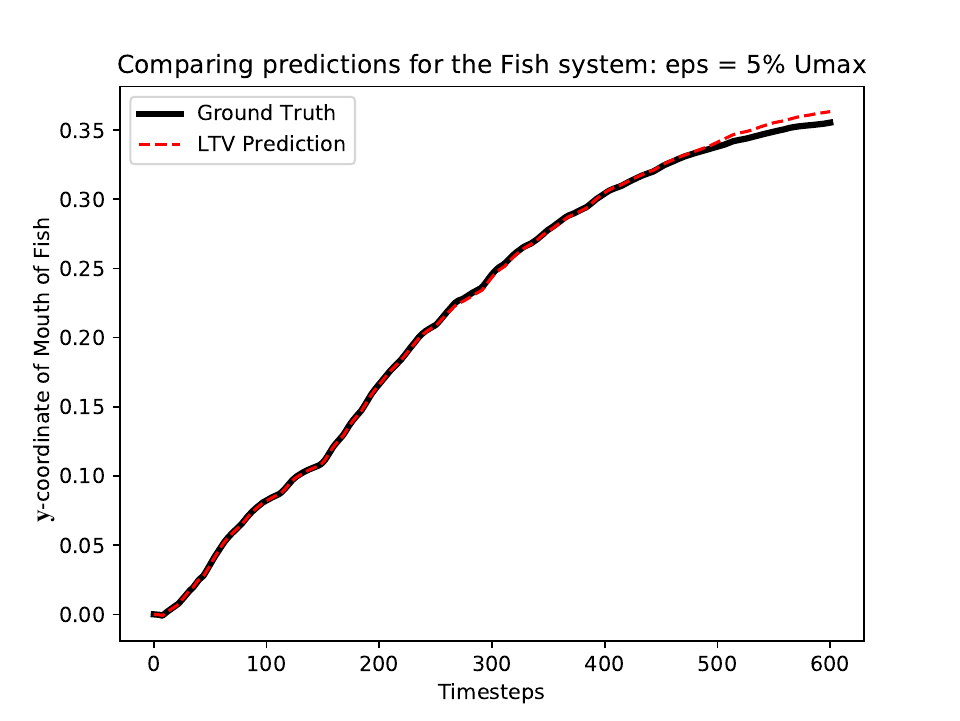}}    
      \subfloat[$30\%$ Perturbation]{\includegraphics[width=1\linewidth]{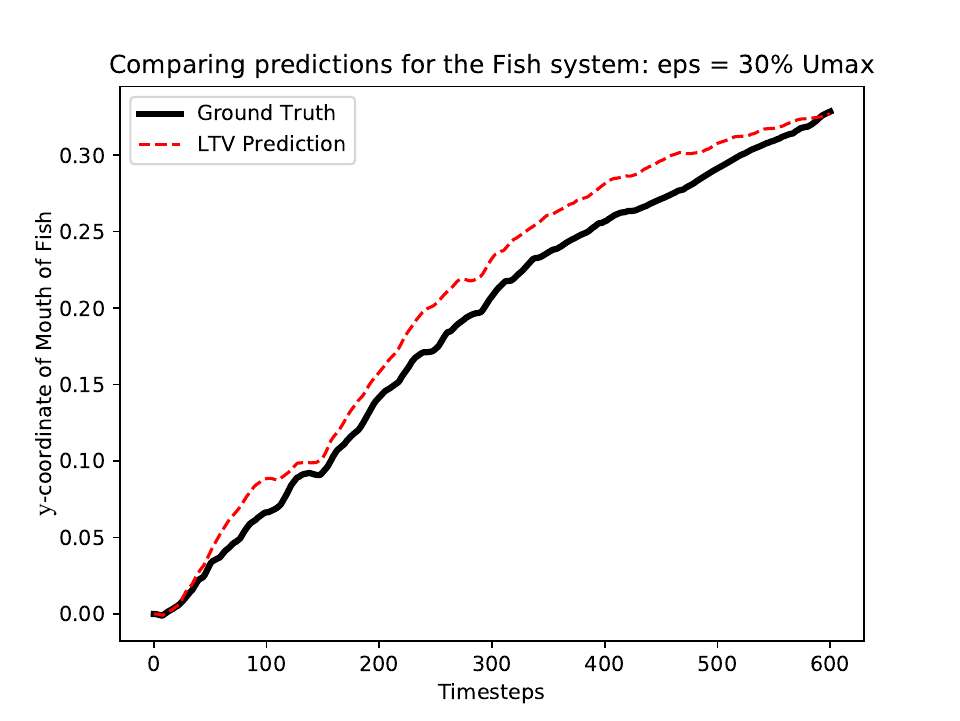}}
      
\end{multicols}
\caption{{Comparing the performance of the surrogate models with the LTV model on Fish trajectories.}}
\label{fig:Fish}
\end{figure}

For the Fish robot, we define the optimal control task of moving from the origin to the target co-ordinates $(0, 0.4,0.2)$ within 600 time-steps. However, the nonlinear surrogate models exhibit significant deviations from the physical system, making them impractical (Fig.~\ref{fig:Fish}). 
In contrast, the recursive local LTV identification coupled with the D2C algorithm successfully drives the Fish robot to the target (Fig.~\ref{fig:optcontrol_fish}). Therefore, we will focus our analysis on the Pendulum and Cartpole systems, which remain practical while still highlighting the identified issues.

\subsection{Model Convergence}
Under a fixed sampling distribution, linear least-squares estimators (SINDy) converge to the Normal-equation solution as $R \rightarrow \infty$ (Prop.~\ref{prop2}). Thus, with a fixed basis and identical sampling distribution, parameter estimates should be invariant across re-samplings in the asymptotic limit. To test this empirically, we draw multiple datasets from the same underlying distribution by varying the random seed, and evaluate all trained surrogates on the same optimal trajectory.

\begin{figure}[!htpb]

\begin{multicols}{2}
      \subfloat[Cart position]{\includegraphics[width=\linewidth]{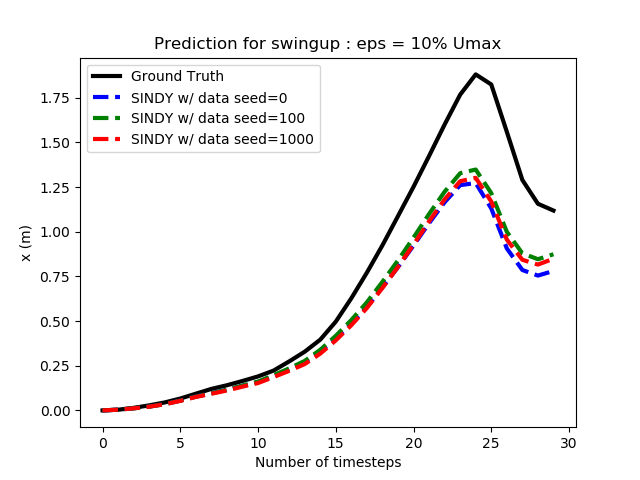}}    
      \subfloat[Pendulum position]{\includegraphics[width=1\linewidth]{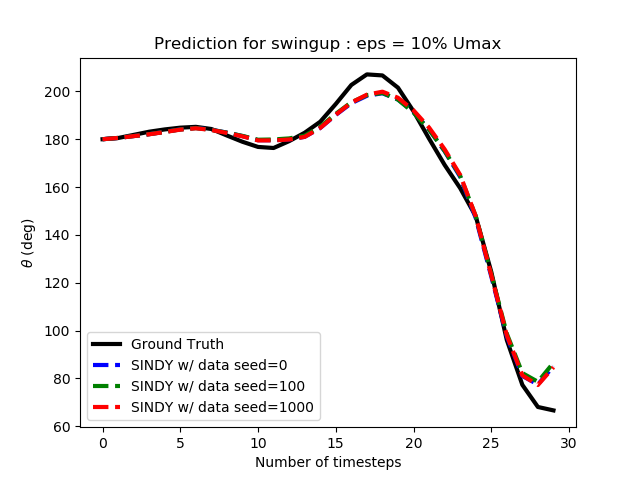}}

\end{multicols}

\begin{multicols}{2}
    \subfloat[Cart position]{\includegraphics[width=\linewidth]{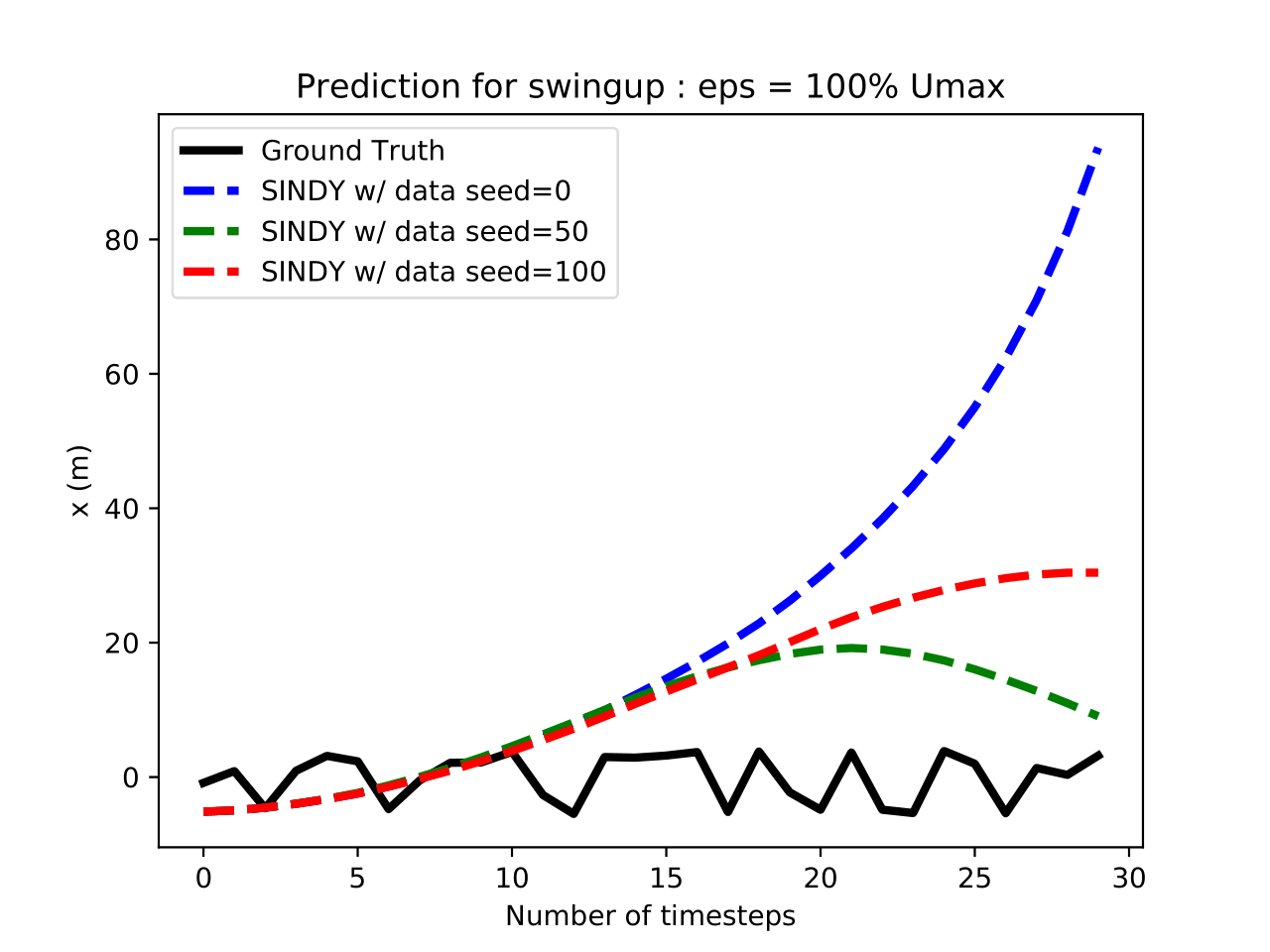}}    
      \subfloat[Pendulum position]{\includegraphics[width=1\linewidth]{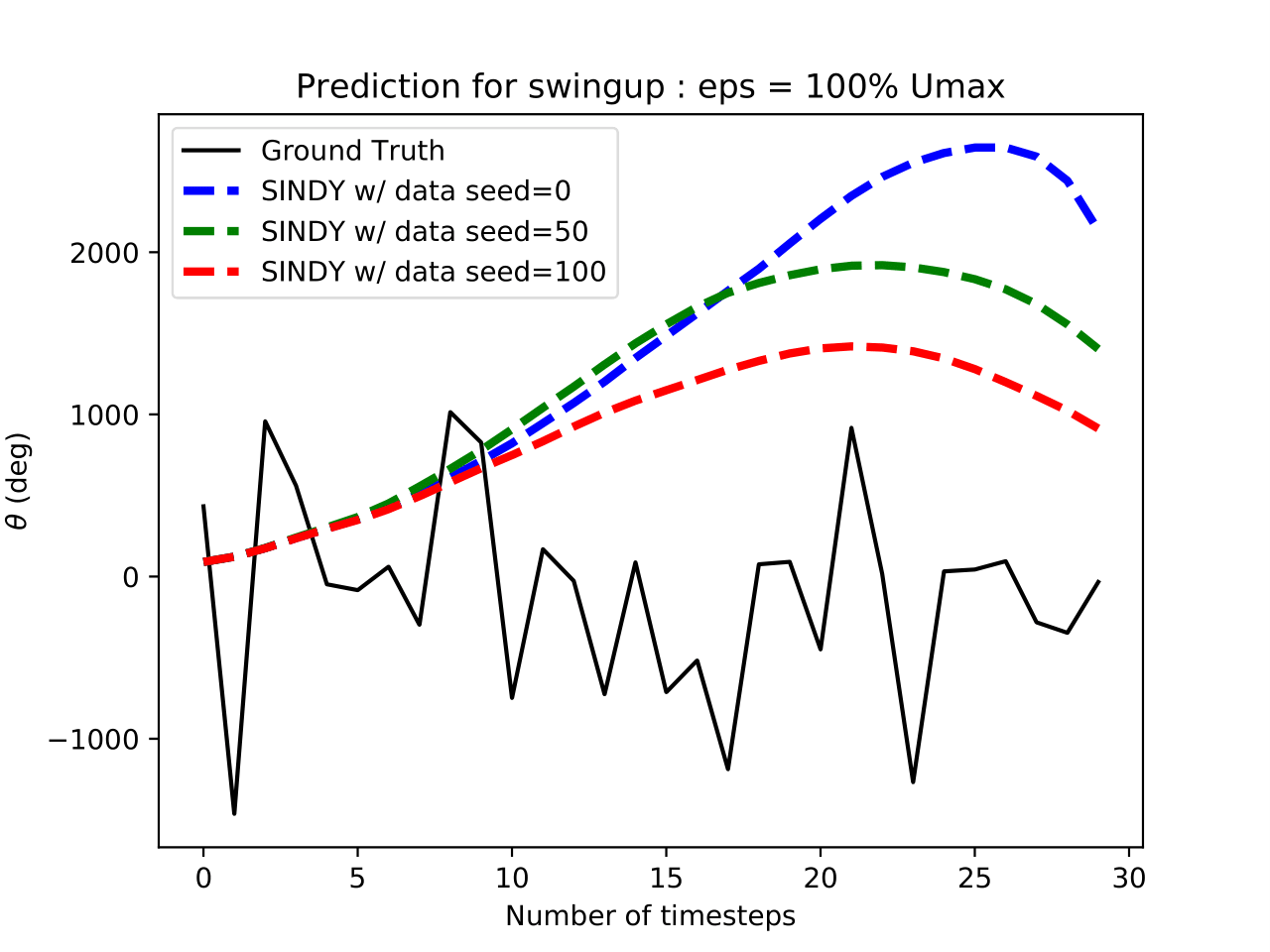}}

\end{multicols}

\caption{{Convergence of the SINDy model learned on datasets with random initial seeds, with perturbations of $10\% U_{max}$ (a)-(b), 
and $100\% U_{max}$ (c)-(d) w.r.t. the optimal trajectory for the Cartpole system.}}
\label{fig:Convergence_SINDY}
\end{figure}

\begin{figure}[!htpb]

\begin{multicols}{2}
      \subfloat[Cart position]{\includegraphics[width=1.1\linewidth]{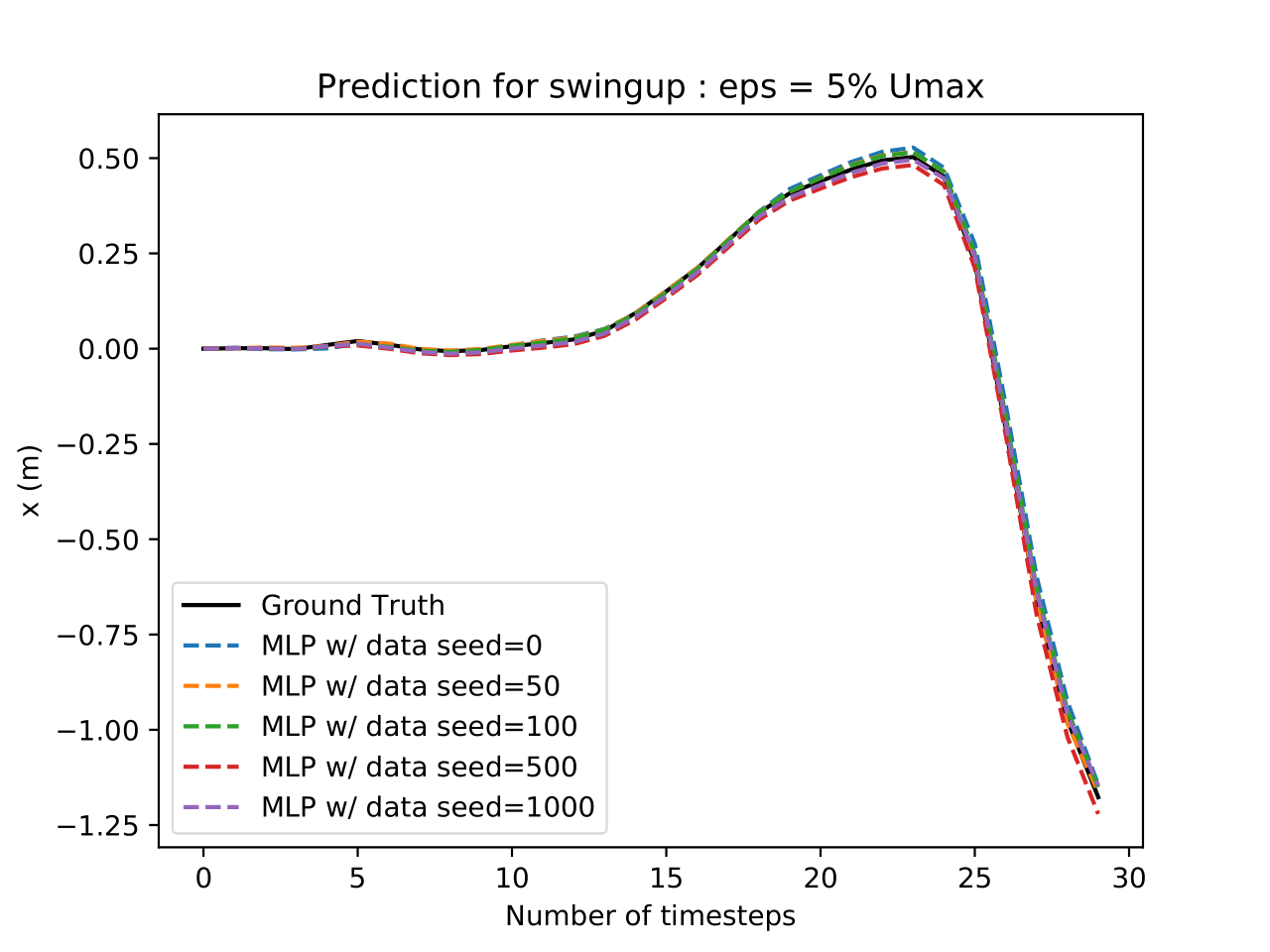}}    
      \subfloat[Pendulum position]{\includegraphics[width=1.1\linewidth]{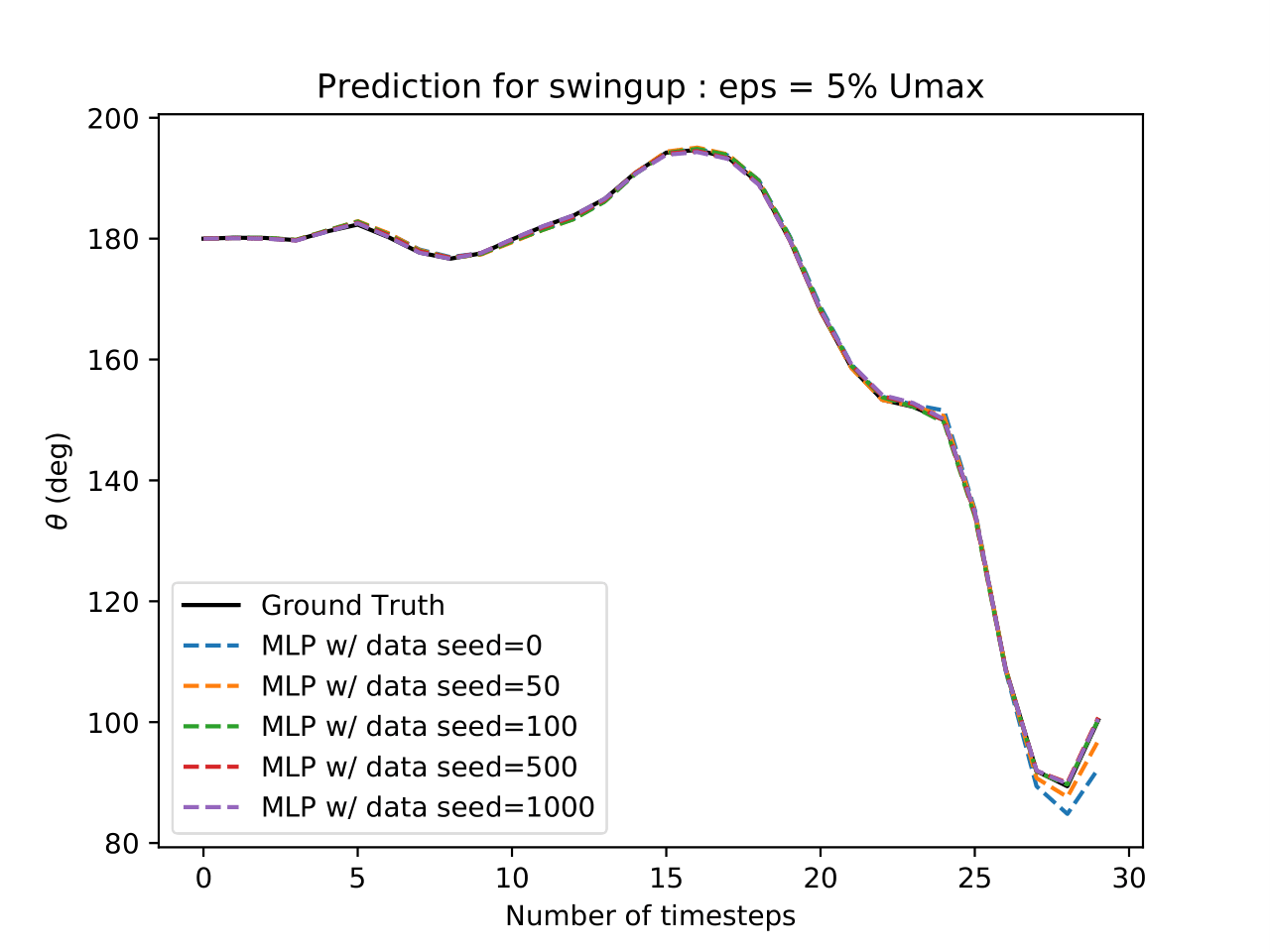}}
      
\end{multicols}

\begin{multicols}{2}
    \subfloat[Cart position]{\includegraphics[width=1.1\linewidth]{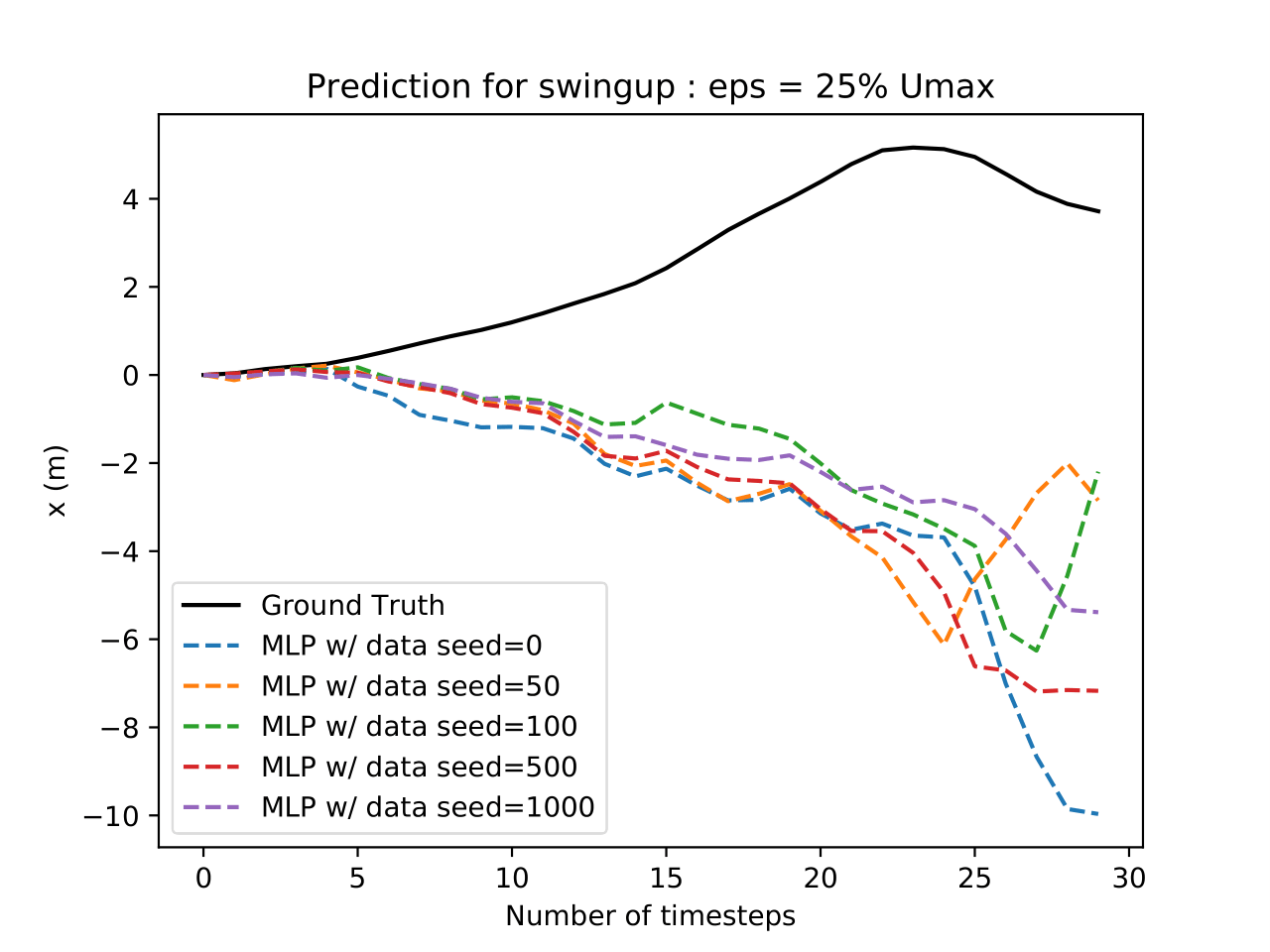}}    
      \subfloat[Pendulum position]{\includegraphics[width=1.1\linewidth]{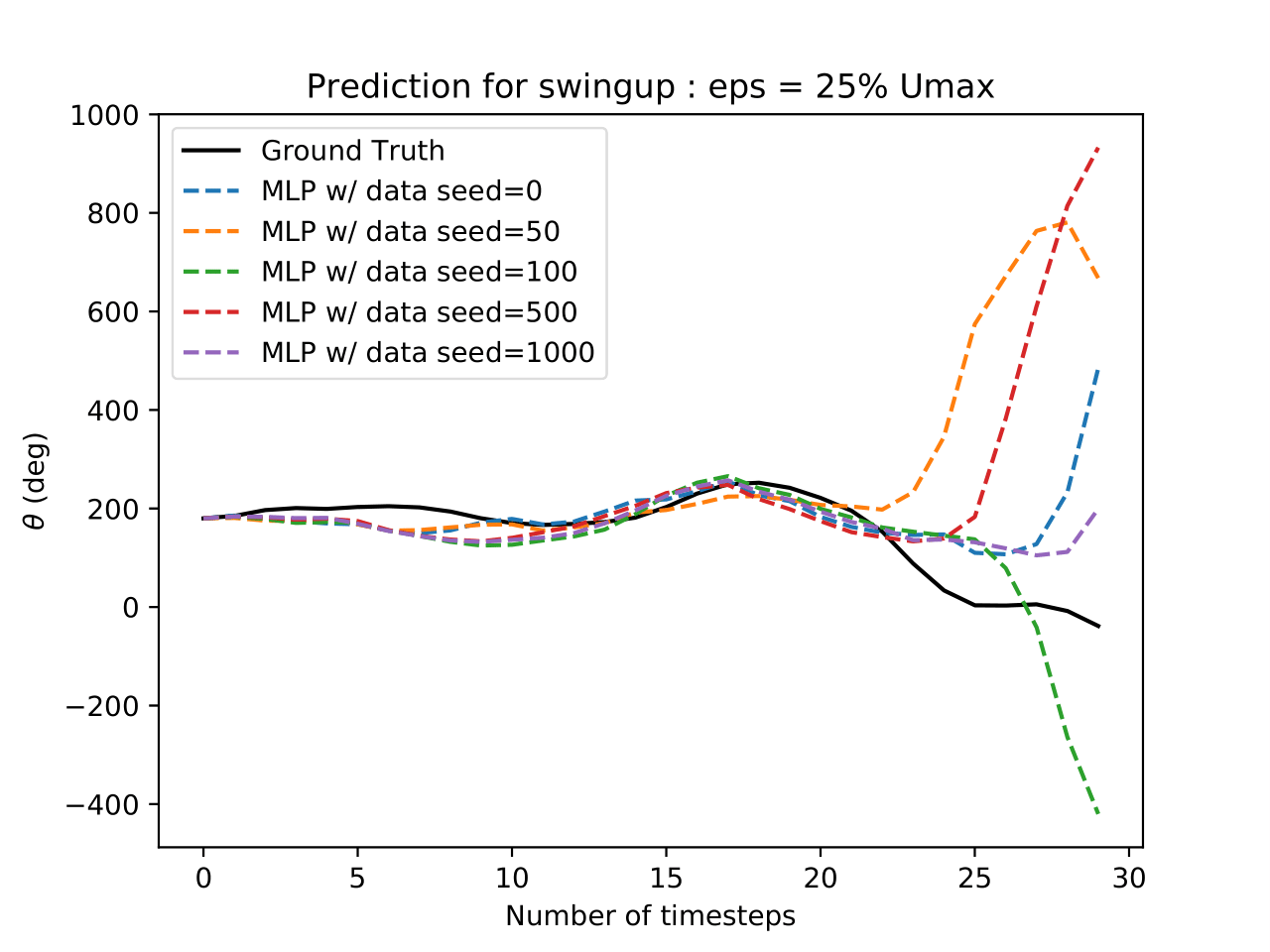}}
      
\end{multicols}

\caption{{{Convergence of the MLP model learned on datasets with random initial seeds, with perturbations of $5\%$ $ U_{max}$ (a)-(b) and $25\%$ $ U_{max}$ (c)-(d) w.r.t. the optimal trajectory for the Cartpole system.}}}

\label{fig:Convergence:MLP}
\end{figure}

In practice, with finite $R$, increasing the sampling neighborhood $\sigma$ causes predictions to diverge across seeds due to the growth of $Var(f \psi^i)$ and the condition number $\kappa(G)$ (Eqs.~\ref{varianceF},\ref{conditioning}). This behavior is evident in Fig.~\ref{fig:Convergence_SINDY}: at higher perturbations, models trained on different seeds yield inconsistent trajectories, unlike the low-perturbation case. For the MLP, which is nonlinear in parameters, convergence to identical models across seeds follows a similar trend (Fig.~\ref{fig:Convergence:MLP}). 

We next quantify the ill-conditioning and variance underlying this lack of convergence, as developed in Sec. \ref{sec_var_illcon}.




\subsection{Ill-Conditioning and Variance}
\begin{figure}[!htpb]
\centering
      \subfloat[Cubic polynomial+Trigonometric basis]{\includegraphics[width=0.8\linewidth]{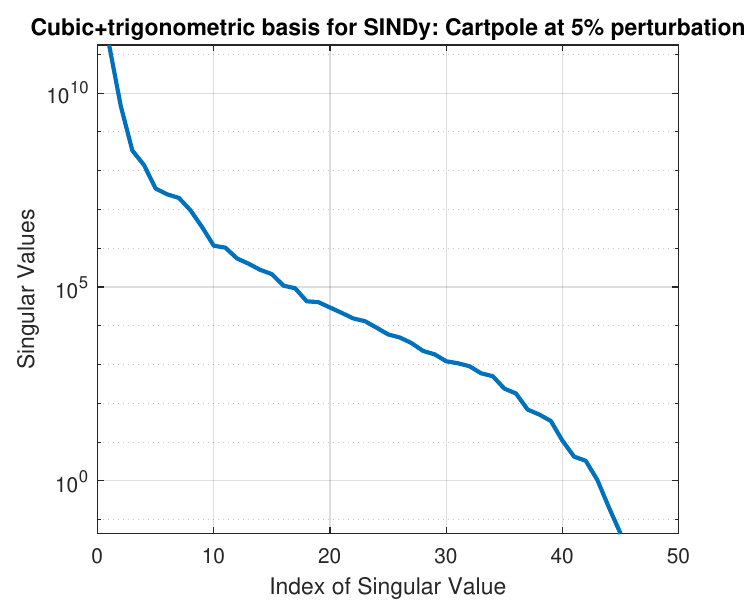}}\\ 
      \subfloat[Fourth order polynomial basis]{\includegraphics[width=0.8\linewidth]{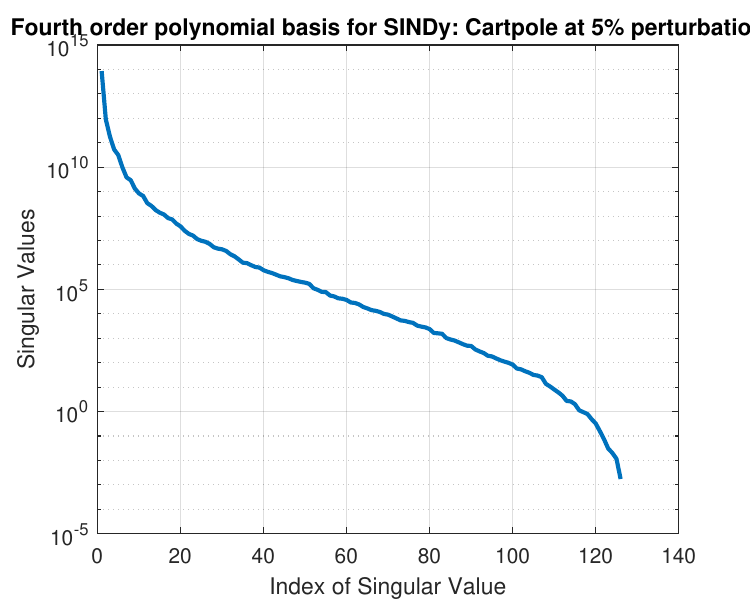}}
      

\caption{{Conditioning of the SINDy model learned on  perturbations of $5\% U_{max}$ w.r.t. the optimal Cartpole swingup trajectory.}} 

\label{fig:IllCond_SINDY}
\end{figure}

For least-squares estimators such as SINDy, we assess ill-conditioning via the empirical Gram matrix $G = (1/R)\Psi^{R'}\Psi^R$, where $\Psi^R_{k,i} = \Psi _i(x^k)$ are basis evaluations on samples $x^k$.
Subsequently, we can utilize the singular-value decomposition to find the conditioning of the Gram matrix (see Fig.~\ref{fig:IllCond_SINDY}). Consistent with the theory (Eq.~\ref{conditioning}), the conditioning number $\kappa(G)$ increases rapidly with the approximation order N (chosen basis) and the exploration level $\sigma$.

\begin{figure}[!htpb]
\centering
      \subfloat[$5\%$ $U_{max}$ perturbations]{\includegraphics[width=0.8\linewidth]{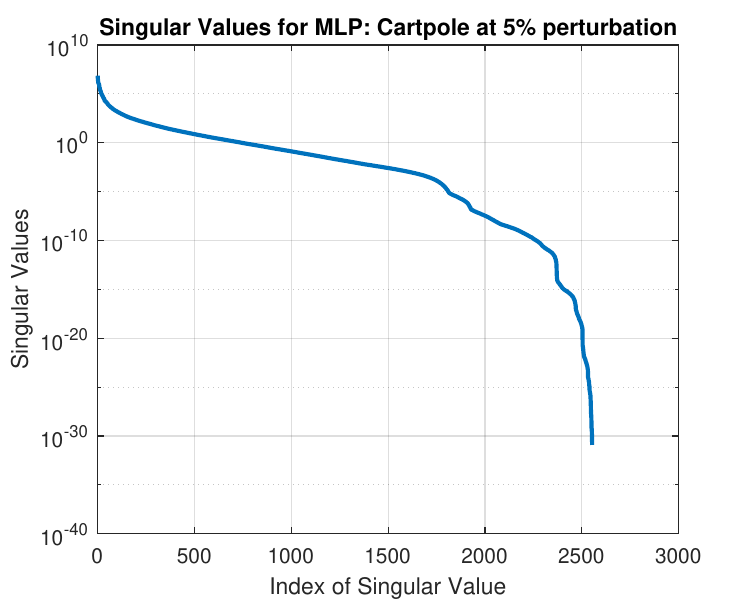}}
      

\caption{{Observing the conditioning of the MLP model learned on  perturbations of 5\% w.r.t. the optimal Cartpole swingup trajectory.}} 
\label{fig:IllCond_MLP}
\end{figure}

Similarly, from Eq.~\ref{eq: LSEquiv}, the conditioning of the Gram matrix $\frac{1}{R}H_t^{R'}H_t^R$ for the MLP algorithm at convergence is studied. We observe that the matrix is highly ill-conditioned in this case as well (Fig.~\ref{fig:IllCond_MLP}), mirroring the linear case and aligning with the $\kappa(G)$ scaling predicted by Eq.~\ref{conditioning}.

\begin{figure}[!htpb]
\centering
      \subfloat[Cubic polynomial basis: $5\%$ perturbation]{\includegraphics[width=0.7\linewidth]{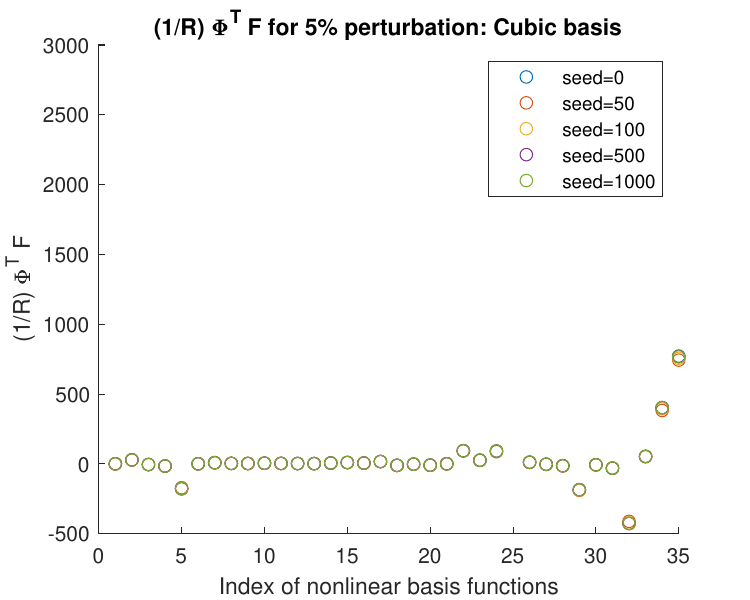}}\\
      \subfloat[Fourth order with trigonometric basis: $15\%$ perturbation]{\includegraphics[width=0.7\linewidth]{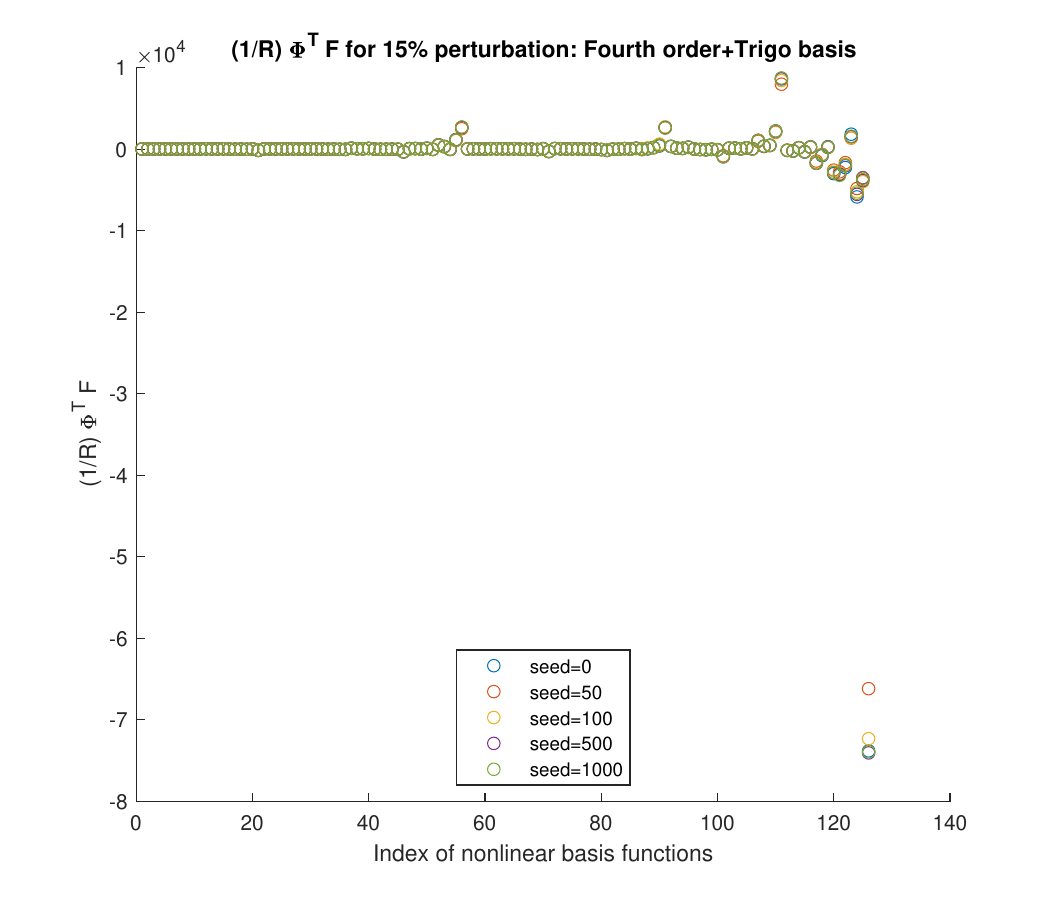}}
      

\caption{{Variance in the SINDy model, learned on  perturbations of $5\%$ and $15\%$ of $U_{max}$ w.r.t. the optimal Cartpole swingup trajectory. 
It is observed that the variance increases with increase in perturbation, as well as increase in order of the basis functions.}} 
\label{fig:Variance_SINDY}
\end{figure}

To check for the variance in estimates from SINDy, we train surrogate models on multiple datasets drawn from the same underlying sampling distribution, and plot the projections of the error w.r.t. the nonlinear basis (see Fig.~\ref{fig:Variance_SINDY}). It is observed that the variance in the estimates scales with both the order of approximation $N$ and the scale of exploration $\sigma$, which agrees with our theoretical development (Eq.~\ref{varianceF}).

\begin{figure}[!htpb]
      \subfloat[At Epoch = 0]{\includegraphics[width=1\linewidth]{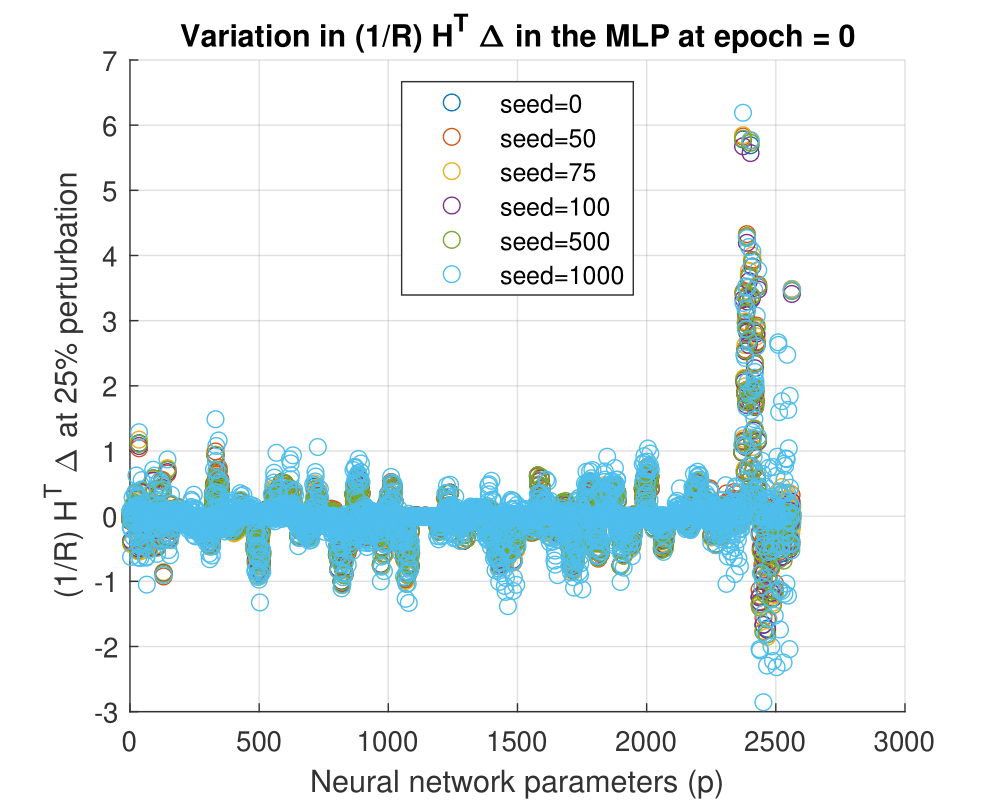}}\\
      \subfloat[At Epoch = 5000]{\includegraphics[width=1\linewidth]{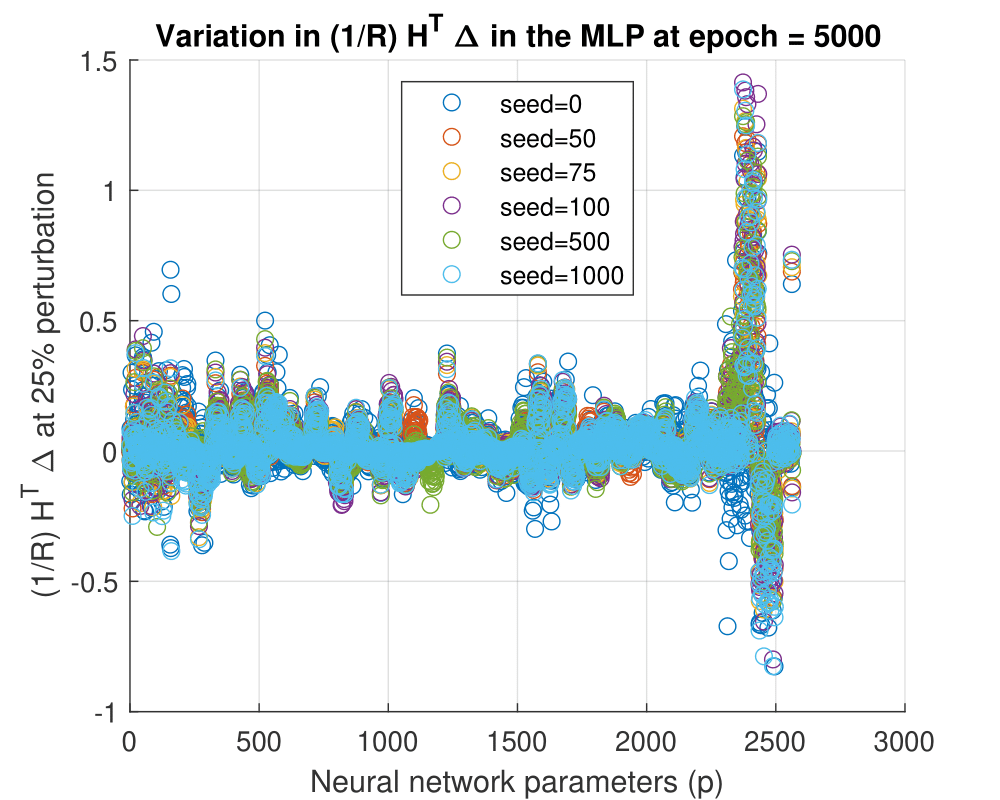}}
      

\caption{{Variance in the MLP model training, learned on  perturbations of $5\%$ of $U_{max}$ w.r.t. the optimal Cartpole swingup trajectory. }}
\label{fig:Variance_MLP}
\end{figure}

Similarly, we train multiple MLP models on each of the data sets, drawn as before from the same sampling distribution, and observe the vector $\frac{1}{R}H_t^{R'}\Delta^R$ corresponding to each of these models, at different points along the iterations leading up to convergence (Fig.~\ref{fig:Variance_MLP}). The plots demonstrate variance at different time steps in the training, even for low to moderate $\sigma$, leading to very different descent directions for different data sets, consequently leading to very different solutions.

\subsection{Optimal Control}

\begin{figure}[!htpb]
      \subfloat[Trajectory evolution]{\includegraphics[width=0.5\linewidth]{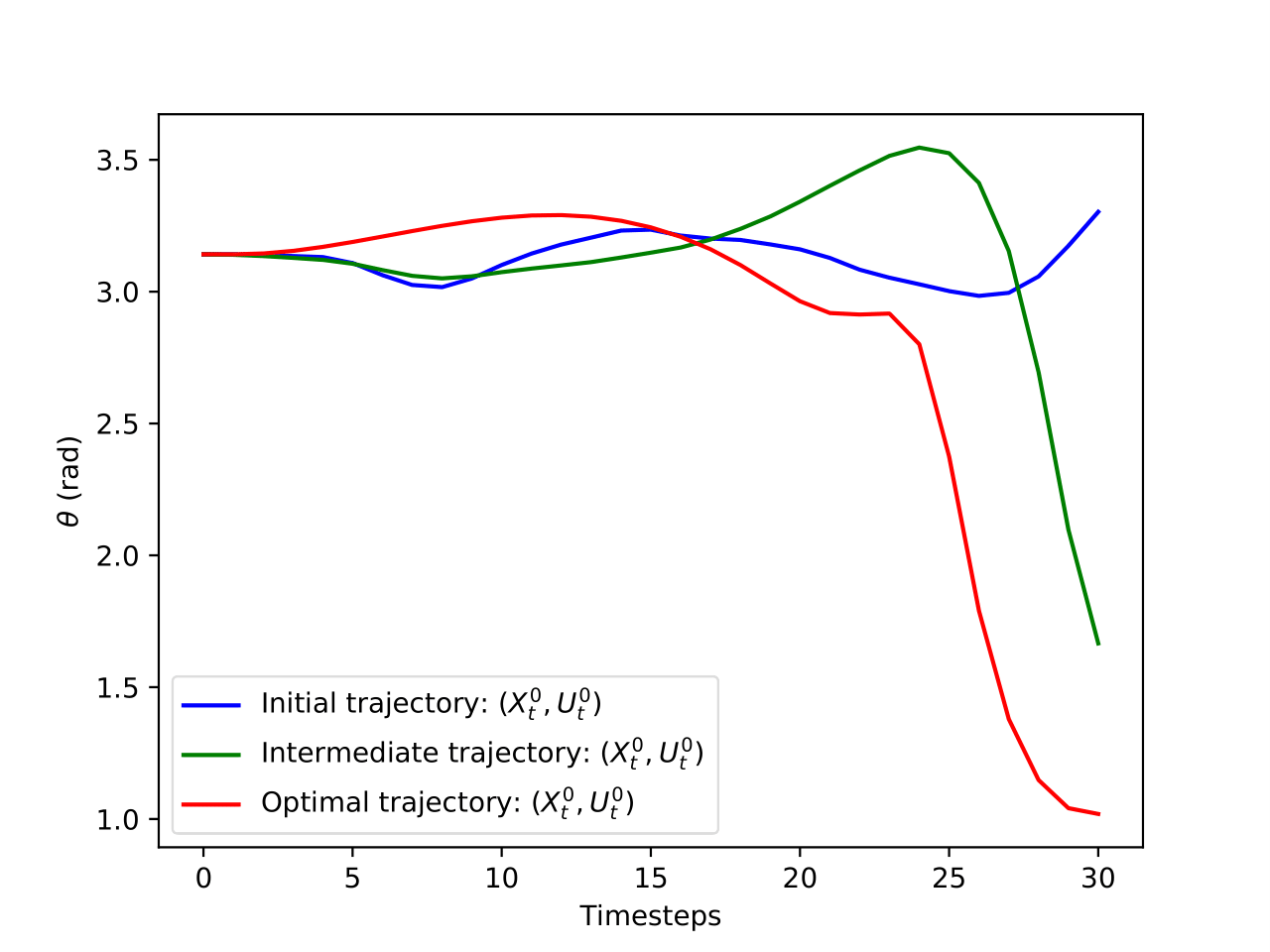}}
      \subfloat[Phase evolution]{\includegraphics[width=0.5\linewidth]{figure/Optimal_ctrl_pendulum_phase-1.png}}\\
      \subfloat[Optimal Control Trajectories]{\includegraphics[width=0.5\linewidth]{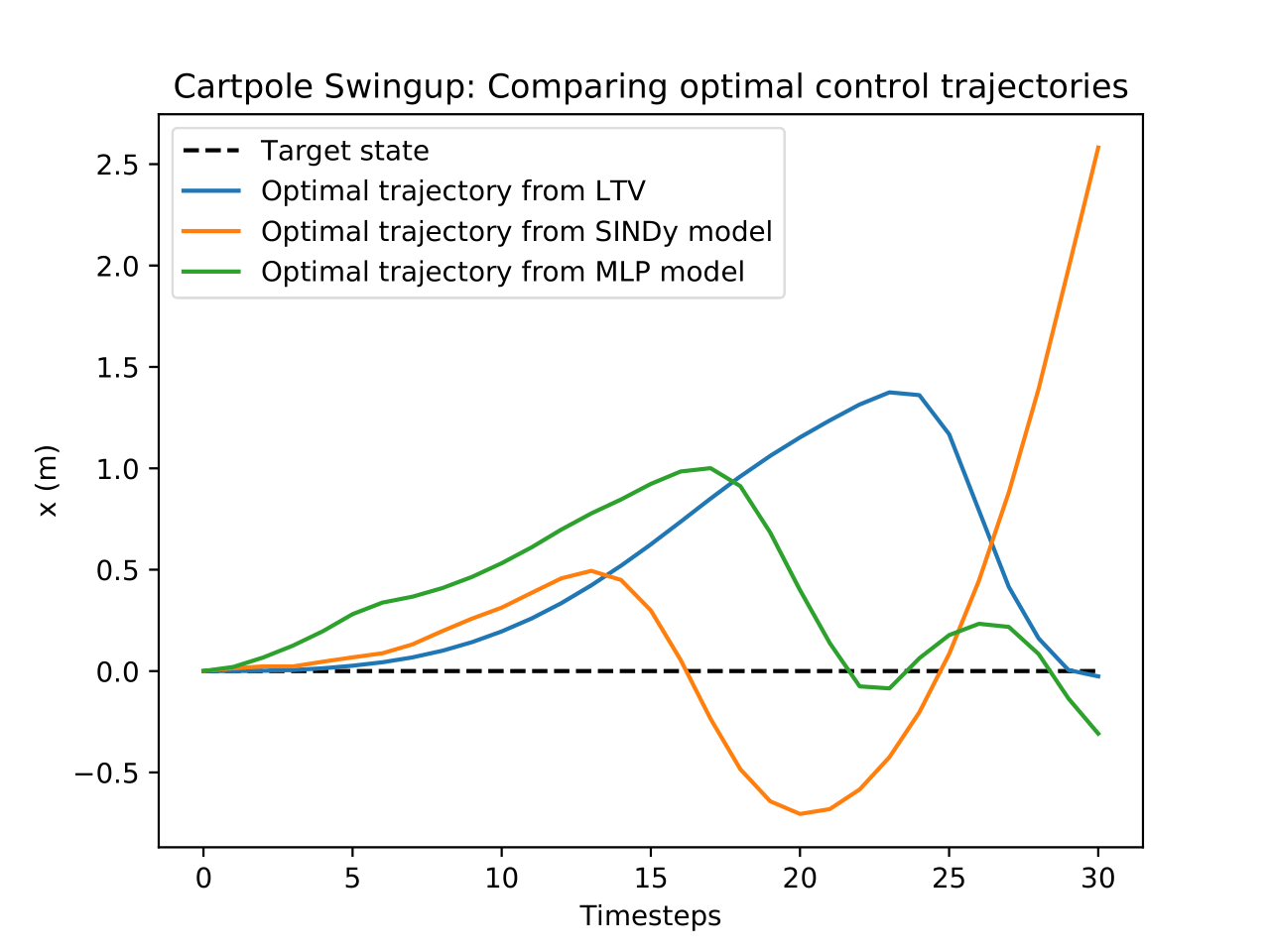}}
      \subfloat[Optimal Control Trajectories]{\includegraphics[width=0.5\linewidth]{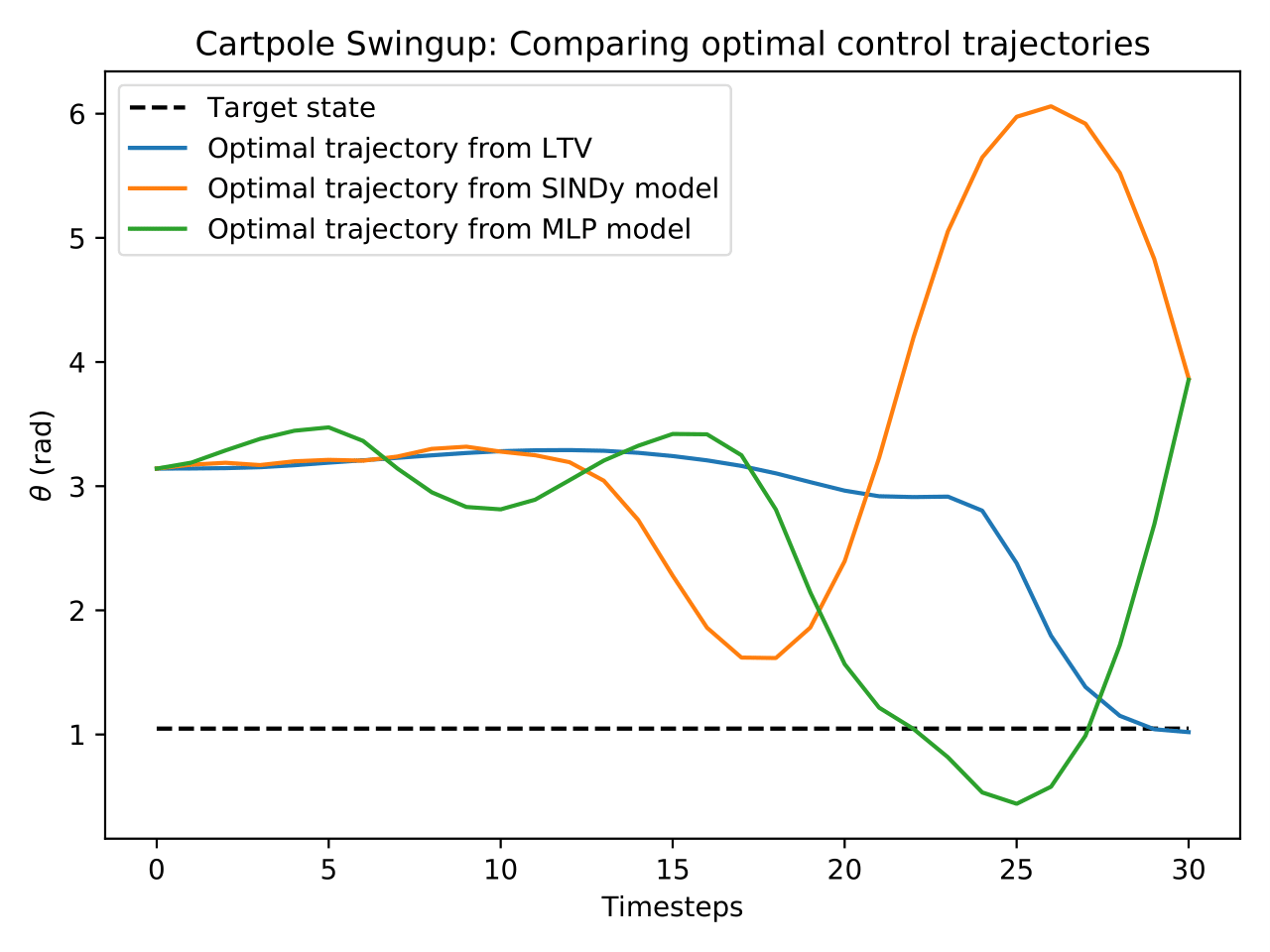}}
      

\caption{{Optimal control of the Cartpole system when using the LTV model as compared to the MLP and SINDy models. The ILQR algorithm, when using the local LTV model can start from a random trajectory and still arrive at the correct solution (a)-(b). However, when using the MLP and SINDy based surrogate models, it fails to converge to the optimal (c)-(d).}}
\label{fig:optcontrol}
\end{figure}

\begin{figure}[!htpb]
      \subfloat[Trajectory evolution]{\includegraphics[width=0.5\linewidth]{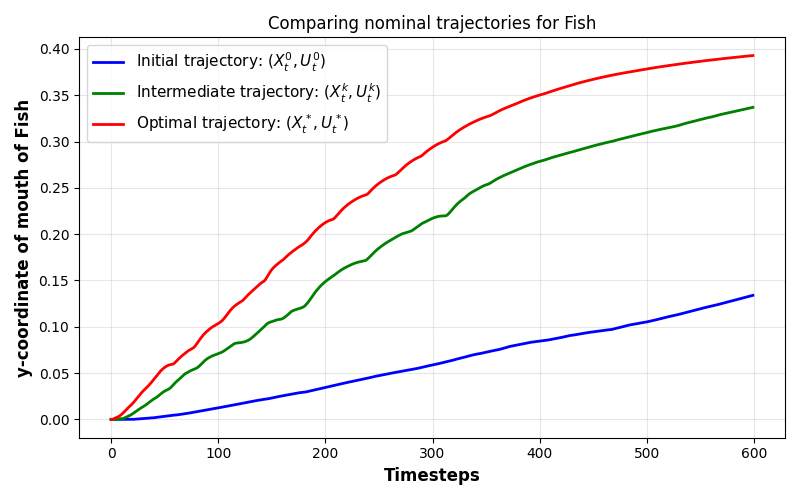}}
      \subfloat[Phase evolution]{\includegraphics[width=0.5\linewidth]{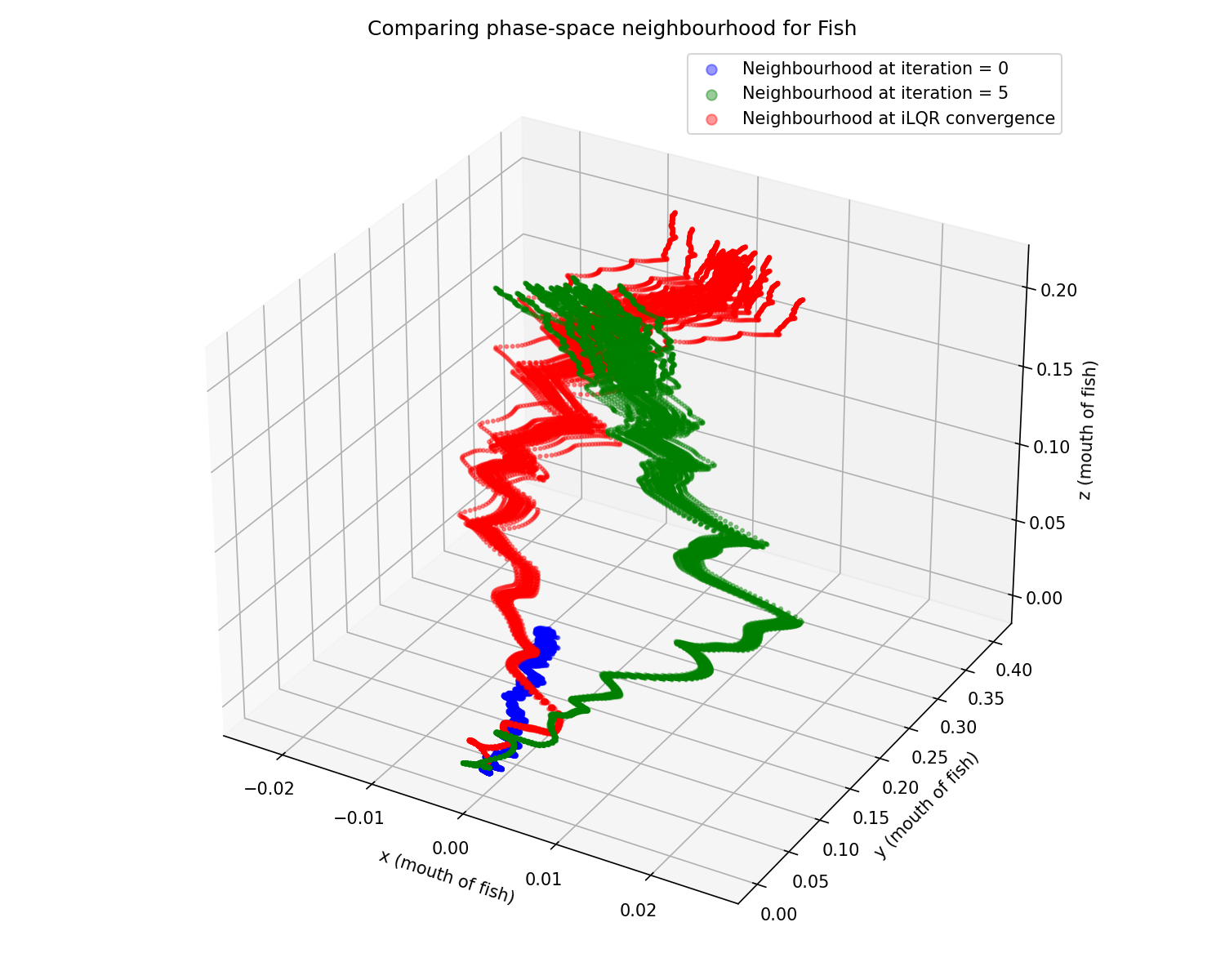}}
      

\caption{{Optimal control of the Fish robot when using the ILQR algorithm with the recursive local LTV system identification. The algorithm successfully drives the Fish robot to the target, thus validating the proposed approach.}}
\label{fig:optcontrol_fish}
\end{figure}

Given that our problem formulation is rooted in solving the optimal control problem, we can now utilize our trained nonlinear surrogate models to compute the optimal control for our corresponding problems, and benchmark them with the approach using a linear time-varying system identification. In general, we observe that, when using the erroneous models learned by the nonlinear system identification techniques, our control algorithm fails to converge to the optimal solution, despite sampling around the true optimum (Fig.~\ref{fig:optcontrol} (c)-(d)). This is in stark contrast to the LTV approach, wherein we can start the sampling anywhere in the state space, and still arrive at the true solution by progressively resampling the state space as the algorithm progresses (Fig.~\ref{fig:optcontrol} (a)-(b), Fig.~\ref{fig:optcontrol_fish} (a)-(b)). Thus, adopting a local LTV system identification approach has an insurmountable advantage over the nonlinear data-driven system identification techniques in the context of optimal control.  

\begin{remark}
    As we have already mentioned in Section 5, the LTV-ILQR based approach was developed by us in the reference \cite{D2C_ASME} and is called D2C. In this reference, an extensive empirical comparison of the approach was done with respect to both state of the art model-based and model-free RL algorithms and the method was shown to be significantly superior in terms of training efficiency, training reliability, optimality and closed loop performance.
\end{remark}

\section{Conclusions}
In this article, we have shown the equivalence of nonlinear model estimation  to the problem of generalized moment estimation of a given underlying sampling distribution, bringing to focus the inherent problems of ill-conditioning and variance naturally incurred by these methods. The empirical evidence clearly demonstrates the ill-conditioning of the Gram matrix, along with the high variance in the forcing function of the constituent Normal equations. Finally, we have also demonstrated that in the context of optimal nonlinear control, the use of a `local' LTV model suffices to combat the issues inherent to nonlinear surrogate models, achieving the global optimum answer, at least locally. \\

\textbf{Global Models:} It is important to clarify that our critique in this paper is aimed at the pursuit of a single, unified `global' model, which we have shown is fundamentally intractable owing to issues of ill-conditioning and variance. This is separate from model-based approaches that construct a comprehensive dynamics model from an ensemble of multiple local models. By restricting identification to local neighborhoods, these methods can successfully alleviate the issues inherent in a single, high-order global approximation. In fact, we can also extend our LTV approximations to do a piecewise ensemble to get a `globally' accurate optimal feedback model. We must note, however, that any attempt to use such a piecewise method to cover the entire state space would quickly become intractable, as it would inevitably encounter the curse of dimensionality. Furthermore, any attempt to interpolate between such models, using say a Neural Network, would also incur the ill-conditioning and variance issues discussed here, in addition to the curse of dimensionality.

\textbf{Model Complexity and Occam's Razor:} It may be tempting to think that a different/ more complex model structure may alleviate the issues of variance and ill-conditioning. However, our analysis suggests that this is going to be an issue regardless of the particular structure, since the nonlinear parametrization assumed is very general and can fit any possible approximating architecture. In the context of model complexity, the LTV model is the simplest and may be construed as a validation of ``Occam's razor" \cite{PAC-2} that the simplest explanation is most likely, given the empirical evidence supports the explanation. This is also consistent with PAC-learning theory which suggests that simpler hypotheses are efficiently learnable \cite{PAC-2}. It is also important to note that the target problem of optimal control/ RL is critical in this context since the lack of knowledge about the optimum necessitates the search over a large domain and to a high order approximation unless done in a local recursive fashion. 

\textbf{Nonlinear vs LTV Architectures.} The great advantage of a nonlinear NN approximation architecture is that it is not limited to a fixed set of basis functions, can be interpreted as continuously modifying the basis functions as the training progresses and can be trained in a highly parallel/ efficient fashion using backpropagation. However, the key disadvantage is that gradient descent cannot be expected to converge to the correct parameter value unlike in linear architectures. Thus, it is not clear that in the context of dynamical systems, the NN model can be expected to converge to the true model using training data as the number of parameters increase, even theoretically. Further, at least in the context of nonlinear optimal control, the issues of variance and ill-conditioning make the model estimation ill-posed unless done locally. Thus, this begs the question whether nonlinear architectures such as deep NNs are useful in the context of nonlinear optimal control? We also stress here that the LTV architecture is required and LTI models do not suffice in this context.\\

\textbf{Partially Observed Problems.} A final note is due here regarding partially observed problems. In general, partially observed problems are significantly more difficult to solve compared to fully observed problems. Typically, the partially observed problems can be ``approximately" converted into fully observed problems by considering a finite past history of outputs and controls, the so called ``information state", leading to a nonlinear autoregressive exogeneous (NARX) model for the output. The NARX model can be represented as: $z_t = \mathcal{N}(Z_{t-1})$, where $Z_{t-1}$ contains a finite $q$-step past history of outputs and control actions. This, in turn, is identical in form to the problem we considered in this paper, and thus, the COVIC will afflict such problems as well. We also mention that the LTV-ILQR approach can be extended to partially observed problems via the notion of the information state as done recently in our work \cite{PODILQR-TNNLS}.

In future work, we plan to extend the analysis and experiments to model-free RL techniques such as DeepRL and show that the issues of ill conditioning and variance that were identified in this paper for model-based RL will be present in the model-free case, if looking for a global approximation of the feedback law. More fundamentally, we are interested in exploring if the ill-conditioning and variance issues identified here for continuous state and control space problems also apply to general supervised/ unsupervised learning problems in discrete spaces such as those in Large Language Models (LLM). \textit{The bottom line conclusion of this paper is that ``more data is not necessarily good, and in fact, can actually be harmful for performance and reliability", at least in the context of nonlinear optimal control problems.}

%% file: Bibliography.tex
\bibliographystyle{IEEEtran}
\bibliography{IEEEabrv,ICRA_refs,TAC_refs}

